\documentclass[11pt]{article}
\usepackage[numbers]{natbib}
\usepackage{url}
\usepackage{enumerate,stmaryrd,mathtools,amssymb}
\usepackage{todonotes}
\usepackage{tikz}
\usetikzlibrary{calc}
\usetikzlibrary{arrows}
\usetikzlibrary{shapes}					

\newtheorem{theorem}{Theorem}

\newcommand{\R}{\mathbb{R}}
\newcommand{\N}{\mathbb{N}}

\begin{document}


\title{Control of multi-agent systems: results, open problems, and applications}

\author{Benedetto Piccoli\\
Department of Mathematical Sciences\\ 
Rutgers University-Camden, Camden, New Jersey, USA, 08102;\\
\texttt{piccoli@camden.rutgers.edu}}

\maketitle


\begin{abstract}
The purpose of this review paper is to present some recent results on the modeling and control of large systems of agents. 
We focus on particular applications where the agents are capable of independent actions instead of simply reacting to external forces. In the literature, such agents were referred to as autonomous, intelligent, self-propelled, greedy, and others. The main applications we have in mind are social systems (as opinion dynamics), pedestrians' movements (also called crowd dynamics), animal groups, and vehicular traffic. We note that the last three examples include physical constraints, however, the agents are able to inject energy into the system, thus preventing the typical conservation of momentum and energy. Also, the control problems posed by such systems are new and require innovative methods. We illustrate some ideas, developed recently, including the use of \emph{sparse} controls, limiting the total variation of controls, and defining new control problems for measures.
After reviewing various approaches, we discuss some future research directions of potential interest. The latter encompasses both new types of equations, as well as new types of limiting procedures to connect several scales at which a system can be represented. We conclude by illustrating a recent real-life experiment using autonomous vehicles on an open highway to smooth traffic waves. This opens the door to a new era of interventions to control in real-time multi-agent systems and to increase the societal impact of such interventions guided by control research.
\end{abstract}

\section{Introduction}
Control theory had a golden age of development after World War II, particularly because of aerospace applications. The Kalman filter, the Pontryagin Maximum Principle, and the Dynamic Programming Principle represent some of the most known peaks of that era \cite{BP07}. 
The applications of control theory went far beyond the aerospace arena, and indeed beyond engineering domains such as the automotive ones. 
In the last decades, control techniques are becoming ubiquitous in many fields such as biology, medicine, and social sciences to mention a few.
The purpose of this review is to provide examples of applications that involve a large number of independent agents, thus posing control problems of a different nature with respect to traditional ones and requiring new methodologies.\\
We refer to a large system when there are one hundred or more agents in a real-world setting. Examples may include large crowds during a public event, a highway with intense traffic, a flock of birds, a social media with many users, as well as drones in flying formation. 
Control problems in such settings may require going beyond classical techniques. On one side,
even a centralized controller may not be able to control every agent of the system. Mathematically this translates into the capability of acting only on specific subsystems. While classical control theory considers this case, change of variables, feedback linearization, and other techniques may be unfeasible because of technological constraints. 
A general demand is that of \emph{parsimonious} controls, i.e. controls that may act on a small number of agents and use minimal information.
Other limitations may come from limited communication and decentralized computation, thus requiring controls to be sufficiently \emph{simple} and having small variation (in time).\\
Systems of this size have been the object of investigations in various physics domains, such as gas dynamics. Several scales can be considered, such as microscopic, mesoscopic, and macroscopic, as well as different points of view, such as Lagrangian and Eulerian. As explained above, for intelligent and autonomous systems, many assumptions may not hold, such as the conservation of momentum and energy. Moreover, typical assumptions of mean-field approaches, such as the indistinguishability of particles, may not hold.
Despite this, various ideas can be inherited from inert matter investigations and used to study these systems. 

\medskip\noindent
To cast our discussion, we start by identifying specific application domains in Section \ref{sec:app}. These applications include vehicular traffic, crowd dynamics, animal groups, and social systems. We first focus on microscopic models, which are often designed via phenomenological modeling. For traffic, we report also macroscopic models which started to be used as early as the 1930s. For most of these problems, mean-field, kinetic, and other mesoscopic approaches are well-known, together with some macroscopic ones. For this reason, we focus on presenting innovative models, which were specifically introduced to deal with such applications.

\medskip\noindent
In Section \ref{sec:new-model}, we start by describing a fully-coupled ODE-PDE system for moving bottlenecks in traffic. The ODE models the time evolution of the bottleneck (e.g. a bus or an autonomous vehicle), while the PDE, a conservation law, captures the time evolution of bulk traffic. The ODE right-hand side depends on the solution to the PDE, while the PDE satisfies a pointwise flux constraint centered at the bottleneck position. Such problems present various new aspects, such as the formation of nonclassical shocks at the bottleneck location, and required new mathematical and numerical techniques.\\
Then we pass to focus on another typical characteristic of vehicular traffic: the presence of multiple lanes. In such a setting, the state of the system must be described by a mixture of discrete and continuous variables. The former describes the lane on which a vehicle is traveling, while the latter represents the position and velocity of the same. Also for this case, new mathematical tools are necessary to study such models.\\
Then we discuss the choice of state spaces and the features of asymptotic behaviors. To keep the presentation short, we provide only a simple example of a social system. Considering the possible multidimensionality of opinions, and their natural bounds, we illustrate an opinion model on a sphere. 
This naturally leads to nontrivial asymptotic states, such as 
\emph{dancing equilibria}. The same richness is encountered in animal groups, for which clusters, lanes, flocks, and mills have been observed.\\
We conclude the section dealing with pedestrian dynamics. The phenomenology in this case includes lane formation in opposite flows, arches at exits, and traffic-light effect at doors. The use of time-evolving measures allows the modeling of microscopic and macroscopic effects at the same time. In turn, this is necessary to reproduce complex phenomena. The natural setting for equations for time-evolving measures is that of Wasserstein spaces, i.e. space of finite-mass measures endowed with the Wasserstein distance.

\medskip\noindent
In the following Section \ref{sec:control-cost} we focus on new emerging control problems. As explained above, large systems call for parsimonious controls. Along this line, we first discuss the concept of \emph{sparse}
control, i.e. affecting a small number of agents, and show results for the celebrated Cucker-Smale model (originally conceived for language evolution but often time interpreted as a model for flocks of birds). Sparsity is typically promoted using $\ell^1$-type costs.\\
Then we focus our attention on the problem of passing to the limit, as in mean-field approaches, in control problems. To pass to the limit in the number of agents and avoid singularities, one is bound to consider highly regular controls, such as Lipschitz continuous. Unfortunately, solutions to optimal control problems are usually much less regular, often times discontinuous. A possible solution is to split the population into a small number of controlled agents, called leaders, and a large number of uncontrolled ones called followers. In this way, one passes to the limit only in the number of followers obtaining a coupled controlled ODE-PDE system. This approach allows $\Gamma$-convergence results \cite{DalMasoBook}.\\
Circling back to the problem of asymptotic behavior, we proposed a new control problem with the aim of avoiding the formation of a cluster. The latter may represent unwanted and rare events, also referred to as \emph{Black Swans.} The considered model is a generalization of the celebrated Hegselmann-Krause model for opinion formation. The results include characterization of cluster avoidance by sparse controls or ineluctability of the same, depending on the strength of interactions among agents.\\
Then we propose another method to promote simple controls, with a small number of switches (change between two values). It is well known that, even for simple systems, optimal controls may exhibit an infinite number of switching. Such phenomenon was proved to be generic in high dimension. The main idea is to introduce relaxed problems by including the total variation (in time) of the control in the cost function. 
Using this technique, one can prove the convergence of optimal controls and provide an estimate of the cost's convergence rate.\\
Lastly, we go back to time-evolving measures and discuss control problems. A possible approach is to look at the whole measure as the naturally evolving quantity and generalize transport equations to the setting of differential inclusions. An alternative one is to stick to the microscopic representation and consider control problems focused on trajectories at the microscopic level, i.e. of the agents. We also consider problems where the group size may vary, such as in confined environments with entrances and exits. In this case, one has to consider a generalization of the Wasserstein distance to deal with measures with different total masses.

\medskip\noindent
In the last Section \ref{sec:future}, we discuss new developments and future perspectives. First, we deal with measure evolution, introducing the new concept of measure differential equations and some generalizations. We also show recent developments for representing the system at different scales and new limiting procedures. Of particular note, the use of \emph{graph-limits} appears to be a promising tool to be used instead of or in combination with the mean-field limit. Another point of view is that of representing solutions as measures in the space of continuous curves. 
For each of these approaches, a mathematical control theory is in its infancy.\\
Finally, we discuss a recent experiment with autonomous vehicles (briefly AVs) to smooth traffic. The idea of using AVs to dissipate stop-and-go waves got a lot of attention in the last decade. First, we discuss a small experiment, with 22 cars, in a confined setting of a ring road. Then we illustrate the recent experiment (November 2022) on the I-24 in Nashville, TN, using 100 AVs. To date, this has been the largest traffic experiment in the world using AVs on an open highway. The feasibility of such an experiment opens the door to other large-scale interventions, in traffic or other domains, which in turn could be the dawn of a new era for the social impact of multi-agent system research.

\section{Motivating applications}
\label{sec:app}
Multi-agent systems may serve the purpose of modeling several different problems where interacting agents are present. The main interests include
identifying simple interaction rules at the local level, which in turn produce interesting self-organization patterns at the group level. The latter may range from lane formation in crossing flows of pedestrians to the remarkable dances of starlings during murmuration \cite{B08,bellomo20review, MR2974186, MR2165531, Degond2013a, Herty2011, Herty2010}. 
Here we focus on three main examples: vehicular traffic, crowd dynamics, and dynamics of social systems.

\subsection{Vehicular traffic}\label{sec:app-traffic}
Vehicular traffic modeling started around one hundred years ago with pioneering work such as that of Greenshields \cite{greenshields}.
The complexity of the problem is
often times underestimated. In fact,
traffic systems include important components to be considered:
\begin{enumerate}
    \item Infrastructure, usually a road network.
    \item Infrastructure regulating technology, e.g. traffic lights.
    \item Physical agents, i.e. the vehicles traveling on the infrastructure. 
    \item Human agents, i.e. the drivers.
\end{enumerate}
Each of the basic components of traffic stimulated research at various levels among engineering, physics, applied mathematics, and other communities. 
An account of this body of research is out of the scope of this paper, but in this section, we will show modeling at two different scales taking into account multiple components.
We refer the reader to \cite{GHPV23,PT12}
for a discussion of models at other scales.
Later, mixed models will illustrate multiscale approaches, and experimental results will show control techniques at work in real life.

\subsubsection{Microscopic models}
\begin{figure}
    \centering
    \includegraphics[width=\linewidth]{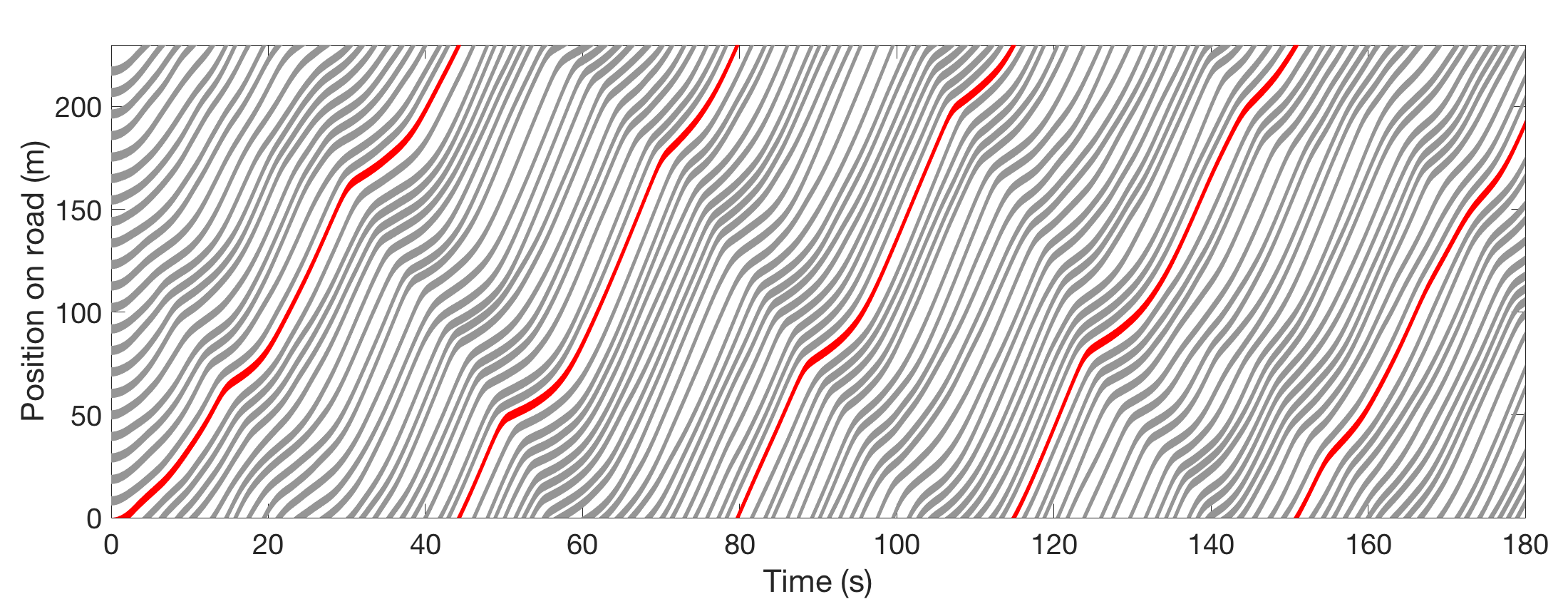}
    \caption{A simulation of Bando-FtL combined model exhibiting appearance and persistence of stop-and-go waves.}
    \label{fig:trajectory_nocontrol}
\end{figure}

Microscopic models are based on set of difference or differential equations for each vehicle in the traffic stream.
Among the most popular models,
we discuss the Follow--the--Leader models  \cite{Gazis1961}, the Intelligent Driver Model (briefly IDM) \cite{Treiber2000}, and the Optimal Velocity or Bando model \cite{Bando1995}.\\
Calling $x_i$, respectively $v_i$, the position, respectively velocity, of the $i$-th driver, the general form of a 
Follow--the--Leader (briefly FtL) model is given by:
 \begin{equation}
\label{model:FtL}
\begin{cases}
\dot{x}_i = v_i,\\
\dot v_{i} = \frac{\phi(v_i,v_{i+1}-v_{i})}{(x_{i+1}-x_{i}-l_{v})^{\gamma}}
\end{cases}
\end{equation}
where $\phi$ is a modeling function,
$l_{v}$ is the vehicle length,
and $\gamma$ is an integer exponent.
The simplest choice is given
by $\phi(x,y)=y$ and $\gamma=1$.
Nice features of FtL models include 
their intuitive formulation and the avoidance of conflicting trajectories.
On the other side, the follower's acceleration will be zero as long as its speed is the same as the leader, but independent from the headway. Thus an
extremely small headway may occur.\\
In the FtL family, the
IDM can be formulated as:
\begin{equation}\label{eq:IDMmodel}
\begin{cases}
\dot{x}_i = v_i,\\
\dot{v}_i = a\bigg(1-\left(\frac{v_i}{v_0}\right)^4 - \left(\frac{s(v_i, v_{i+1} -v_i)}{x_{i+1}-x_{i}-l_{v}} \right)^2\bigg),
\end{cases}    
\end{equation}
where $a$ is the maximal acceleration,
$v_0$ is the maximal speed, and
$s$ is a modeling function. IDM has been widely used in engineering work, but it presents a number of mathematical drawbacks as illustrated in \cite{IDM22}.\\
The Bando model has a different formulation focusing on the attempt of the driver to achieve an ideal velocity, which depends on the headway:
\begin{equation}
\label{model:Bando}
\begin{cases}
\dot{x}_i = v_i,\\
\dot{v}_i = \bigg(V(\Delta x_i)-v_i\bigg),
\end{cases}
\end{equation}
where the modeling function $V(\cdot)$
is called the optimal velocity function
and, in the original model, was chosen as:
\begin{equation}\label{defB}
V(x)=V_{\max}\frac{\tanh(\frac{x-l_{v}}{d_{0}}-2)+\tanh(2)}{1+\tanh(2)},
\end{equation}
where $d_0$ is the safety distance and $V_{\max}$ the maximum velocity.\\
The combination of Bando and FtL models \cite{Gong2022}
is amenable to reproducing the emergence of traffic waves  observed in real traffic 
\cite{GHPV23}.
To illustrate this fact,
we consider 22 vehicles on a ring road with dynamics given by a combination of the FtL model \eqref{model:FtL} and the Bando ones \eqref{model:Bando} with weights
$\alpha =0.5 \text{m}/\text{s}^2$ for FtL and $\beta = 20\text{s}^{-1}$ for Bando. 
To use realistic parameters, the maximal speed is set to $9.75 \mathrm{m/s}$ and the vehicle length to $l_v= 4.5 \mathrm{m}$. The obtained simulation is reported in Figure \ref{fig:trajectory_nocontrol} with one trajectory highlighted in red.
One can observe the appearance and persistence of so-called stop-and-go waves observed in the famous experiment \cite{sugiyama2008} performed in Japan. 

\subsubsection{Macroscopic models}
Arguably the most famous fluid dynamic model is the Lighthill-Whitham-Richards (briefly LWR) one consisting of a single conservation law:
\begin{equation}
\rho_t+f(\rho)_x=0
\end{equation}
where $\rho$ is the car density and $f=\rho\cdot v $ the flux function with $v$ the average speed. The model is obtained by assuming the conservation of the car quantity and $v=v(\rho)$.\\
There is a rich history of designing models with more conserved quantities,
including the Payne-Whitham model
\cite{Payne1971,Whitham74}, famous Daganzo's criticism
\cite{DaganzoCritic}, and the Aw-Rascle model \cite{AwRascle},
which, in conservation form, can be formulated as:
\begin{equation}
\left\{
\begin{array}{c}
\rho_t + \partial_x (\rho\,v)=0      \\
y_t + \partial_x (y\, v)=0
\end{array}
\right.
\end{equation}
where $y=\rho (v+p(\rho))$ is a modified momentum and $p$ an appropriate pressure term. The model avoided the drawbacks of the Payne-Whitham one. It was extensively studied, and various generalizations were provided, see \cite{PT12}. \\
Starting from the seminal paper by Holden and Risebro \cite{HR95}, many macroscopic models were extended to networks, represented  by topological graphs. Fairly complete accounts of results are available, see 
\cite{GPbook,GHPbook}.

\subsection{Crowd dynamics}
The term crowd is used in different ways in the literature, including specifically referring to pedestrians, or more generally to groups of agents having \emph{social} behavior. The latter includes decision-making processes taking into account an energy functional, as well as other criteria dictated by the interaction with other agents.
To unify the perspective, here we first present this section with some models for pedestrians, then turn attention to animal groups.

\subsubsection{Pedestrian models}
Pedestrian movements were modeled somehow later than other systems, but starting in the early 2000s a lot of models started to appear.
We recall the popular
\emph{social force} one proposed
by Helbing and collaborators
\cite{Helbing1995} inspired by
the work of Lewin \cite{Lewin1951}. 
This is a microscopic model:
\begin{equation}\label{e-HM}
\left\{
\begin{array}{l}
\dot{x}_i=v_i\\
\dot{v}_i=F_i(v_i)+\sum_j F_{ij}(x_i,x_j,v_i,v_j)+F_{E}(x_i,v_i),
\end{array}
\right.
\end{equation}
where $x_i$ is the position of the $i$-th agent, $v_i$ its velocity, $F_i$ a term similar to Optimal Velocity for traffic,
$F_{ij}$ captures interactions among agents, and $F_E$ the interaction with the environment. 
As for traffic, different explicit expressions have been proposed for the so-called force terms $F_i$, $F_{ij}$, and $F_E$, and we try to summarize the main features.
The term $F_i$ drives the pedestrian to a desired velocity, which may be composed of a desired direction, depending on the final destination, but also a desired speed, which depends both on the pedestrian characteristic (e.g. age) and the environment.  
The terms $F_{ij}$ mainly model trajectory conflict resolution, but may also contain attractive terms.
Finally, the term $F_E$ is capturing the need to avoid obstacles and  move within a structured environment.
The terms  $F_{ij}$ and $F_E$ 
can be associated to potentials $\Phi$ and $\Psi$ by setting:
\begin{equation}
F_{ij}(x_i,x_j,v_i,v_j) =\nabla \Phi(x_j-x_i, v_j-v_i), \quad
F_{E}(x_i,v_i)=\nabla\Psi(x_i,v_i).
\end{equation}
This model has been extensively used, 
see \cite{Bellomo2016,CPTbook}, both to capture the phenomenology of pedestrian dynamics and to describe normal behavior versus panic. The latter is defined as a different behavior observed during emergencies, see 
\cite{Cercignani1994,Cristiani2014,Maury2011}. We provide some examples of pedestrian formations in Section \ref{sec:meas-evol}
The model was also improved and corrected using experimental data, see
\cite{Arechavaleta2008,Chitour2012,Farina2016}.

\subsubsection{Animal group models}
\begin{figure}[h!]\centering
\includegraphics[width=4.2cm,angle=90]{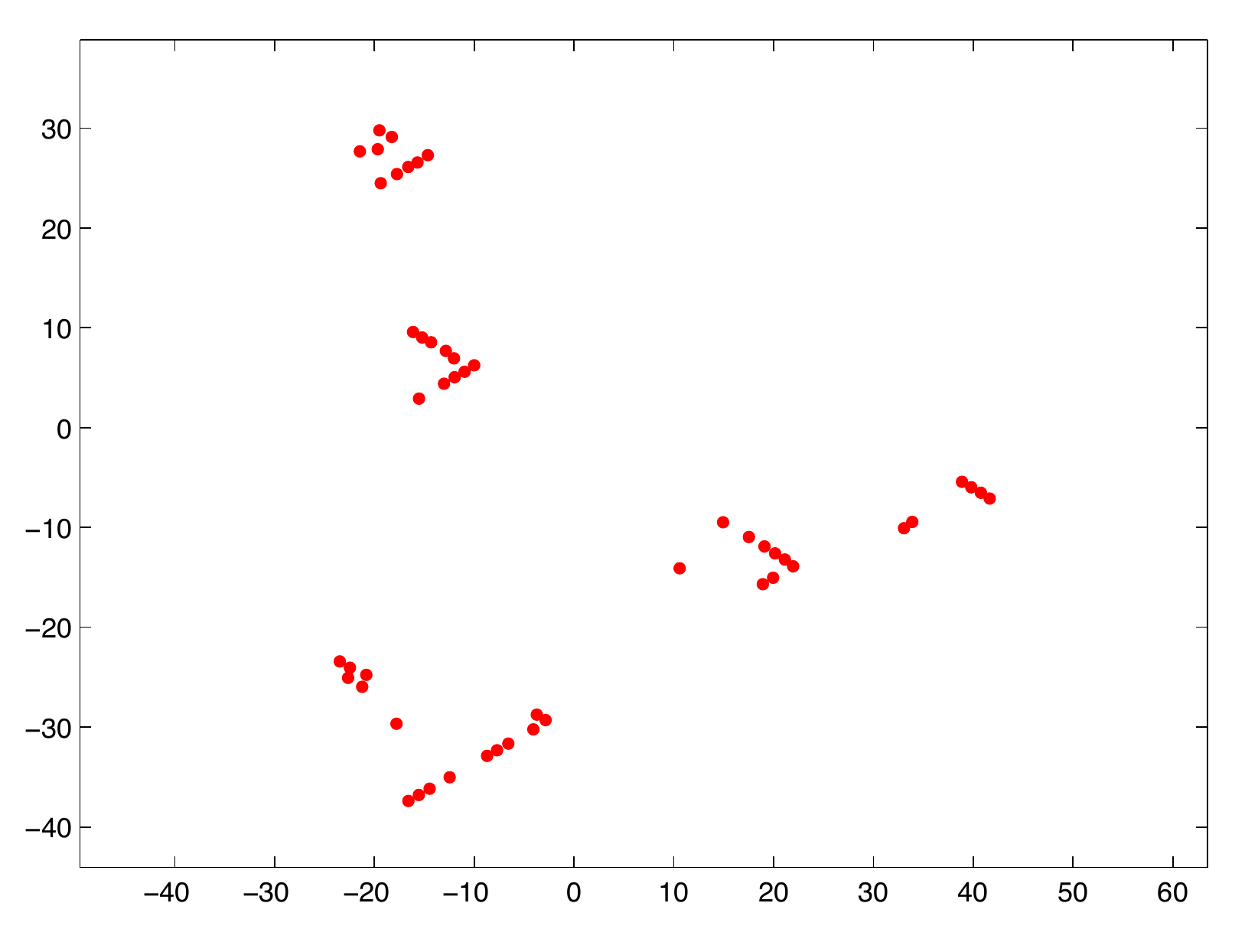}
\includegraphics[width=6cm]{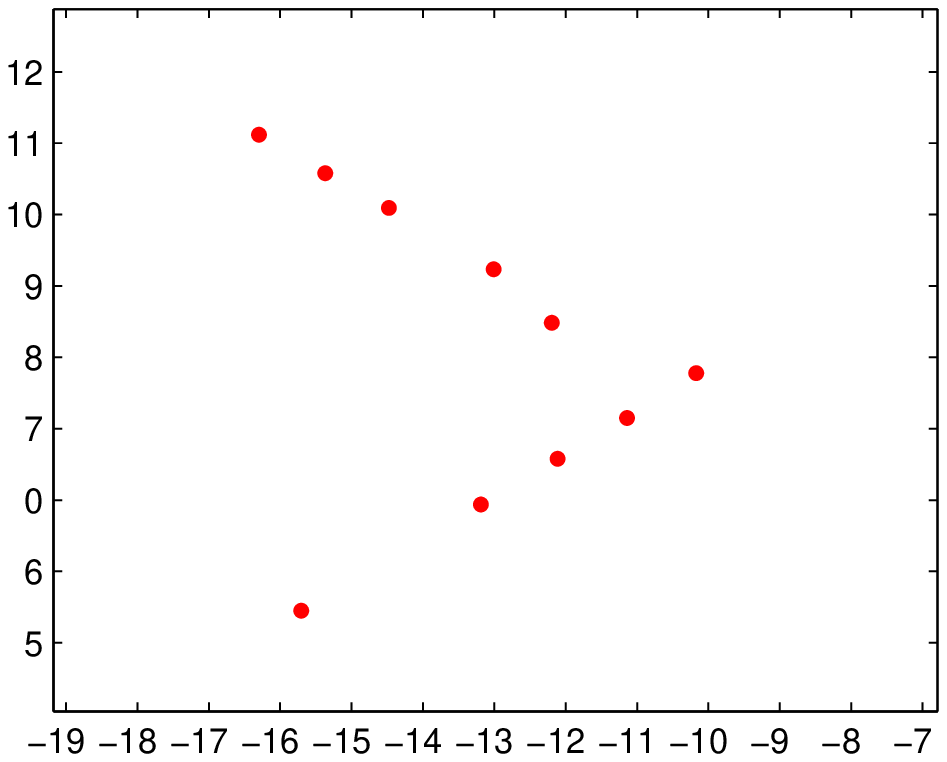}
  \caption{V-like and echelon formations obtained with model \eqref{eq:animal-gen}: large group traveling left to right (left), zoom on a V-like formation (right).}
  \label{fig:Vee}
\end{figure}

Several models were proposed in the specialized literature on animal groups. Indicating again by $x_i$ the position of the $i$-th agent,
a general model takes the form:
\begin{equation}\label{eq:animal-gen}
\dot{x}_i= \sum_{j\in A_i} \phi_a(x_i,x_j) \frac{x_j-x_i}{\|x_j-x_i\|}+
\sum_{j\in R_i} \phi_r(x_i,x_j) \frac{x_i-x_j}{\|x_j-x_i\|}
\end{equation}
where $A_i$ is the set of agents interacting with the $i$-th agent by attraction, $\phi_a$ is a modeling function depending on agent positions, and similarly for $R_i$, $\phi_r$, and repulsion. Notice that this model is first order since the velocity is prescribed directly. The choice is justified by the fact that, despite being called social forces, the agent is able to modify its speed almost instantaneously.
The functions $\phi_a$ and $\phi_r$ may contain terms that depend on the angle between agents $x_i$ and $x_j$ (in two or three dimensions) and are related to the visual field of the animal under consideration.
Such models are able to reproduce various self-organization patterns observed in nature, such as lines, crystal clusters, V shapes, and others, see Figure \ref{fig:Vee}.
For example, elephants form lines, as well as lobsters on the ocean floor, while ducks form crystal clusters on rivers and cranes V shapes when migrating.
We refer the reader to \cite{CFP11} and references therein for a more extensive discussion.

\begin{figure}[t!]
    \centering
    \includegraphics[width=\textwidth]{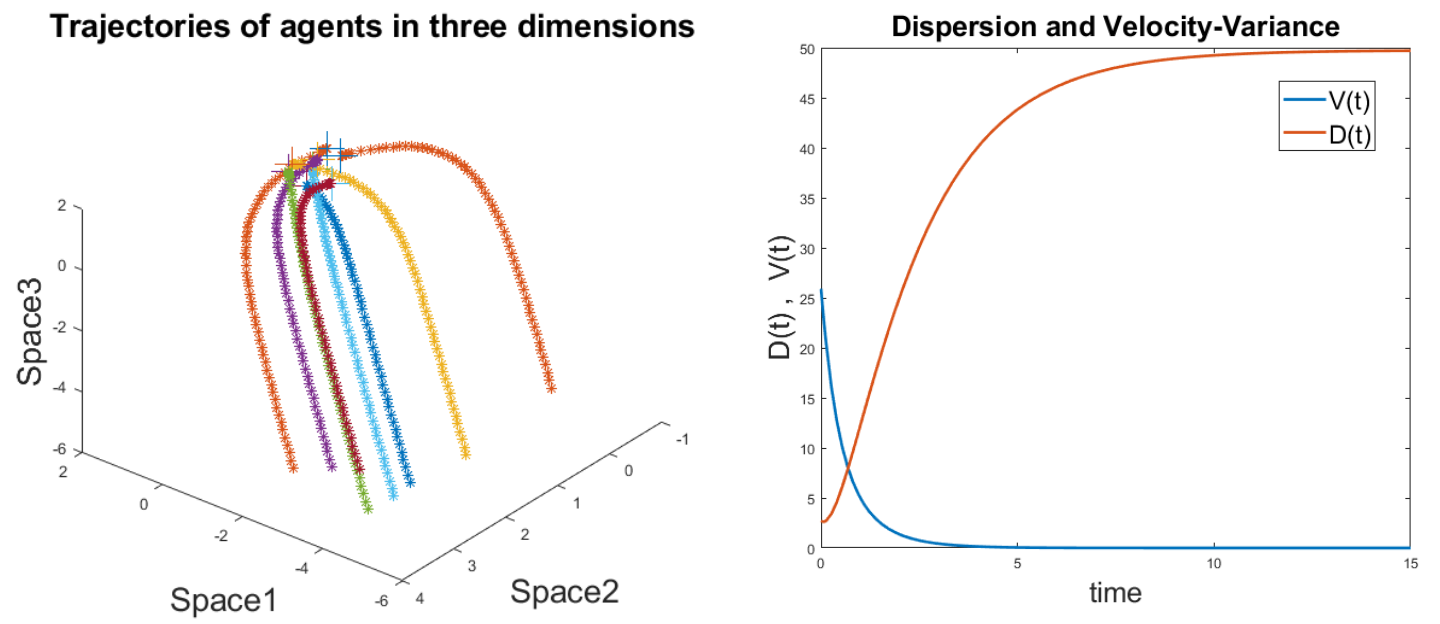}
    \caption{Left: time evolution of the Cucker-Smale model~\eqref{eq:CS} with $a(s)=\frac{1}{(1+s^2)^{\frac13}}$ in $\R^3$.
    Right: time evolution of $D=\sum_{ij}\|x_i-x_j\|^2$ and
    $V=\sum_{ij}\|v_i-v_j\|^2$ measuring position and velocity dispersion of the group.}
\label{fig:CS1}
\end{figure}

While many models were considered by the biology community, a specific model attracted the attention of many applied mathematicians and physicists:
the \emph{Cucker-Smale model} (briefly CS) \cite{Cucker2007}. Interestingly,
the model was originally designed for the evolution of languages, but it reproduces a specific dynamic feature, called \emph{alignment}, which is also typical of birds. The model for $N$ agents with positions $x_i$ and velocities $v_i$ reads:   
\begin{equation}\label{eq:CS}
\begin{cases}
\dot{x}_i(t)&=v_i(t)\\
\dot{v}_i(t)&=\frac{1} {N}\sum\limits_{j=1}^N \phi(\Vert x_{j}(t) - x_{i}(t)\Vert )(v_j(t)-v_i(t)),
\end{cases}
\end{equation}
where $\phi$ is modeling the dependence of attraction on the distance among agents. 
Usually, $\phi$ is decreasing as the distance among agents increases, as in the original choice $\phi(r)=\frac{1}{(1+r^2)^\beta}$ with $\beta>0$.
For $\beta<\frac{1}{2}$, the agents asymptotically align their velocity, see Figure \ref{fig:CS1}. 

\subsection{Social dynamics}
\label{sec:app-soc}
Social dynamics is a fairly general term, which may include crowd dynamics.
Here we are interested in presenting a simple example of opinion dynamics.
The celebrated Hegselmann-Krause model (briefly HK, \cite{Hegselmann2002})
models the evolution of opinions, represented by a single parameter,
in a group of agents interacting only with proximal neighbors.
This characteristic is called \emph{bounded confidence} and mathematically corresponds to restricting interactions of the $i$-th agents only with those within a given radius, e.g. on the set
$\{y: \|y-x_i\|\leq 1\}$. This results in a dynamic that asymptotically converges to a finite number of clusters at distance greater or equal to $1$.
In formula:
\begin{equation}\label{eq:HK}
\dot{x}_i =  \sum_{j: \|x_i-x_j\|\leq 1}(x_j-x_i).
\end{equation}
The set of interacting agents changes with the evolution, thus the problem is indeed nonlinear and even discontinuous in time and space.
Given the discontinuity in spaces, one has to carefully choose the appropriate concept of solution. The most common approach is to use Caratheodory solutions, but this implies non-uniqueness, and statements about the asymptotic states have to be formulated meticulously. We refer the reader to \cite{CFPR21} for an extensive discussion and examples.
The HK model gave rise to many variants and generalizations, see
\cite{Jabin2014,Motsch2014,PR21}.


\section{New modeling approaches}
\label{sec:new-model}

%
%
%

%

In this section, we present models which were introduced in the last ten-fifteen years
and presented new features w.r.t. those presented in Section \ref{sec:app}.
The first one, motivated by traffic, exhibits a full coupling between an ODE, representing a moving bottleneck, and a PDE, for the bulk traffic evolution.
The second, still motivated by traffic, is a hybrid model with continuous variables, position, and speed of vehicles, as well as discrete ones, the lanes on which vehicles travel.
The third one, motivated by social dynamics, is defined on a nontrivial manifold (the $d$-dimensional sphere) and produces interesting asymptotics.
Finally, we describe how equations for time-evolving measures were implied in pedestrian dynamics.

\subsection{ODE-PDE models for moving bottlenecks}
\begin{figure}[h]
  \centering
    \begin{tikzpicture}[line cap=round,line join=round,x=1.cm,y=0.7cm]

      \draw[<->] (0.2, 6.) -- (0.2, 1.) -- (5.3, 1.);

      \draw (0.2, 1.) to [out = 80, in = 180] (2.35, 4.5)
      to [out = 0, in = 100] (4.5, 1.);

      \draw (0.2, 1.) to [out = 80, in = 180] (1.77, 3.5)
      to [out = 0, in = 100] (3.5, 1.);

      \draw[dashed, very thick] (0.2, 1.) -- (0.2 + 4.5, 1. + 1.7);
      \draw[dashed, very thick] (0.2, 3) -- (0.2 + 4.5, 3. + 1.7);

      \draw[dashed] (4.2, 2.5) -- (4.2, 1.);
      \draw[dashed] (3.25, 4.17) -- (3.25, 1.);
      \draw[dashed] (0.78, 3.2) -- (0.78, 1.);

      \draw[dashed] (1.52, 3.4) -- (1.52, 1.);

      \node[anchor=north,inner sep=0] at (0.2,0.9) {$0$};
      \node[anchor=east,inner sep=0] at (0.2,3) {\small $F_\alpha(\omega)$};
      \node[anchor=north,inner sep=0] at (4.1,0.9) {$u_\omega$};
      \node[anchor=north,inner sep=0] at (1.52,0.9) {$\tilde{u}_\omega$};
      \node[anchor=north,inner sep=0] at (3.25,0.9) {$\hat{u}_\omega$};
      \node[anchor=north,inner sep=0] at (0.78,0.9) {$\check{u}_\omega$};
      \node[anchor=north,inner sep=0] at (5.3,0.9) {$\rho$};
      \node[anchor=east,inner sep=0] at (0.2,6) {$f$};
      \node[anchor=north west,inner sep=0] at (2.2,3) {$f_{\alpha}$};

      \node[anchor=south east,inner sep=0] at (4.8,4.8) {$\varphi_\omega$};
      \node[anchor=south east,inner sep=0] at (4.8,2.8) {$\omega \rho$};

    \end{tikzpicture}    
  \caption{Graphical representation of $F_\alpha$, $\tilde{u}_\omega$, $\check u_\omega$, and $\hat u_\omega$.}
  \label{fig5:check-hat}
\end{figure}
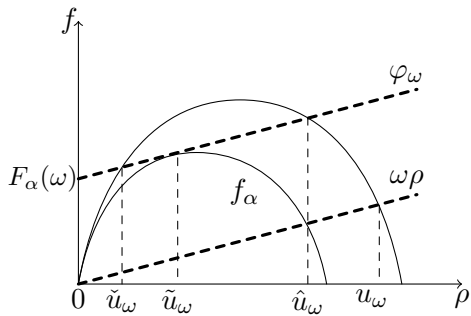

In vehicular traffic, a moving bottleneck can model a large bus or truck or
a vehicle driving differently from an autonomous one.
Some mixed-scale models coupling ODEs and PDEs were proposed, see \cite{DMG2014exist,Lattanzio2011},
together with numerical schemes \cite{Chalons2017,DMG2014num,Gasser2013}. 
and various extensions \cite{GG2011, PiacentiniGoatinFerrara2019platoon}.
We focus on the approach proposed in~\cite{DMG2014exist, GGLP2019, LLG98}
given by:
\begin{equation}\label{eq:DMG-control}
\left\{
\begin{array}{l}
\rho_t + f(\rho)_x=0     \\
\dot y(t) = \min\left\{ \omega(t), v(\rho(t,y(t)+))\right\}\\
f (\rho(t,y(t))) - \dot y(t) \rho (t,y(t)) \leq 
F_\alpha\left(\dot y(t)\right)\\
\end{array}
\right.
\end{equation}
where $f(\rho)=\rho \,v(\rho)$,
the control $\omega\in [0,U]$ indicates a speed chosen by the bottleneck, say an autonomous vehicle or a vehicle with autonomous cruise control, 
$v(\rho(t,y(t)+))=\lim_{h\to 0, h\geq 0}
v(\rho(t,y(t)+h))$,
and $F_\alpha$ is given by:
\begin{equation} \label{eq5:constraintdef}
F_\alpha\left(\dot y(t)\right):=\max_{\rho} \left( \alpha f(\rho/\alpha)-\rho\dot y (t)\right),\qquad \alpha \in\ ]0, 1[.
\end{equation}
The function $F_\alpha$ 
measures the capacity drop given by the presence of the  autonomous vehicle driving at a different speed, usually slower. 
To compute the value of the maximum in the definition of $F_\alpha$,
we first consider the concave function $f(\rho/\alpha)= \frac{\rho}{ \alpha } v(\frac{\rho}{\alpha})$, see
Figure \ref{fig5:check-hat}.
For every $\omega$ let $\tilde{u}_\omega$ be such that
$f'_{\alpha}(\rho) = \omega$,
 and set
 $\phi_\omega(\rho)=
f_{\alpha}(\tilde{u}_\omega) + \omega \left(\rho - \tilde{u}_\omega\right)$. 
If $\dot{y}=\omega$, i.e. the autonomous vehicle drives slower than the traffic average speed, then 
$F_\alpha = \phi_\omega(0)= 
f_{\alpha}(\tilde{u}_\omega) - \omega \tilde{u}_\omega$. Otherwise,
if $\dot y(t)=v(\rho(t,y(t)+))$
then the last inequality 
in \eqref{eq:DMG-control} is obvious
since the right-hand side vanishes.
This model produces nonclassical shocks corresponding to values 
$\hat{u}_\omega$ and $\check{u}_\omega$
defined as follows.
Set $I_\omega = \left\{\rho :\, f(\rho) = \phi_\omega (\rho)\right\}$, then
$\check{u}_\omega=\min I_\omega$
and $\hat{u}_\omega=\max I_\omega$.
We refer again to Figure \ref{fig5:check-hat} and to \cite{BDMGPbook,GGLP2019}
for details. Moreover, continuous dependence results are available for such a system
\cite{LiardPiccoli19,LiardPiccoli21}.

\subsection{Hybrid model for multilane traffic}
Many of the classical models for traffic, as in Section \ref{sec:app-traffic},
deal with single-lane roads or averages over different lanes. However, lane-changing occurs on highways, and other multilane roads, and is a big component of traffic phenomena. For instance, lane-changing is suspected to be a trigger of traffic waves, and is a cause of capacity drop effect \cite{YKLH17}.\\
Lane-changing maneuvers are usually either voluntary or mandatory, see \cite{laval2006lane,jin2013multi,ZHENG201416} for details. We focus on the former type of lane change and consider two conditions for its occurrence:
\begin{description}
    \item[Safety:] There should be enough space in the new lane to allow the maneuver, also considering velocities of vehicles;
    \item[Incentive:] The new lane is less crowded or, in any case, has better traffic conditions.
\end{description}
There are other possible factors to consider, including visibility, time from the last lane change, the position of the lane, etc.\\
Following \cite{MOBIL07}, we focus on modeling the two conditions via constraints on acceleration. 
For a fixed lane $k$ and vehicle $i$,
we indicate by $a^k_i(t)$ its acceleration at time $t$, and by  $\bar{a}_i^{k'}(t)$ 
the expected acceleration if the lane is changed to $k'$ (assuming $k$ and $k'$ are adjacent). The conditions can be stated as:
\begin{description}
   \item[incentive condition:] $\bar{a}_i^{k'}(t) \geq a_i^k(t) + \Delta$
   \item[safety condition:] $\bar{a}_i^{k'}(t) \geq -\Delta \ \text{ and } \ \bar{a}_{j(i,k')}^{k'}(t) \geq -\Delta.$
\end{description}
where  $\Delta$ is a positive modeling parameter and $j(i,k')$ is the vehicle immediately behind the expected new position of vehicle $i$ on lane $k'$.
In simple words, the first condition requires that the acceleration in the new lane $k'$ exceeds the gain $\Delta$.
The second requires the acceleration of vehicle $i$ and the one immediately behind in lane $k'$ to be not too negative, i.e. both vehicles should not need to brake excessively.\\
The obtained model requires nontrivial notation  to be detailed. In particular, a complete state of the system requires the use of mixed variables $(k,x_k,v_k)(t)$ for each vehicle, where $k$ is the lane and $(x_k,v_k)$ the position and velocity on lane $k$.
Notice that to specify a Follow-the-Leader-type dynamics, one has to reconstruct for each vehicle  
the position of the leader in the current lane. Finally, we obtain a hybrid stochastic system and
we refer the reader to \cite{GPV21} for details.

\subsection{State space and asymptotic states in social dynamics}
\label{sec:sphere}
\begin{figure}[h!]
    \centering
\includegraphics{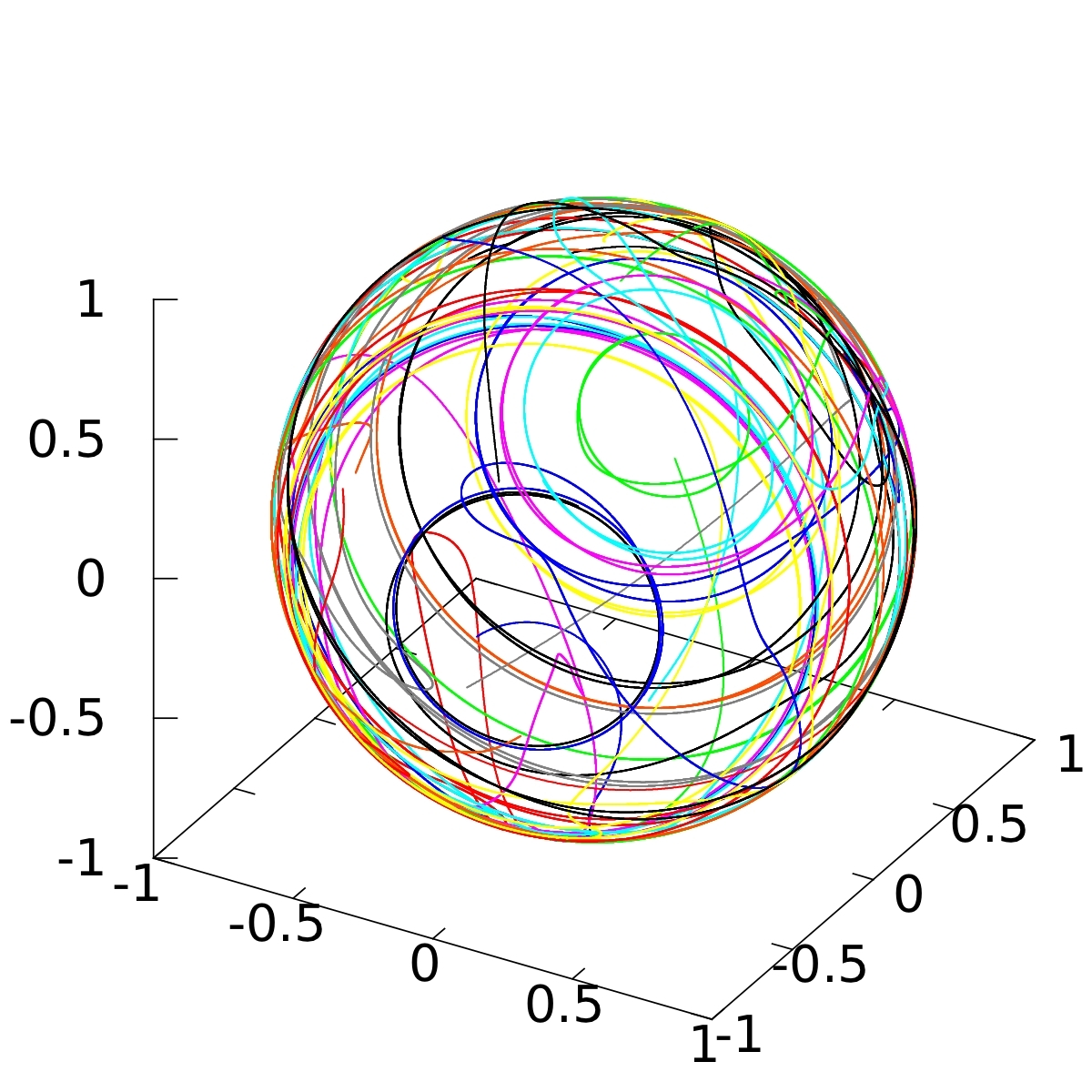}
\includegraphics{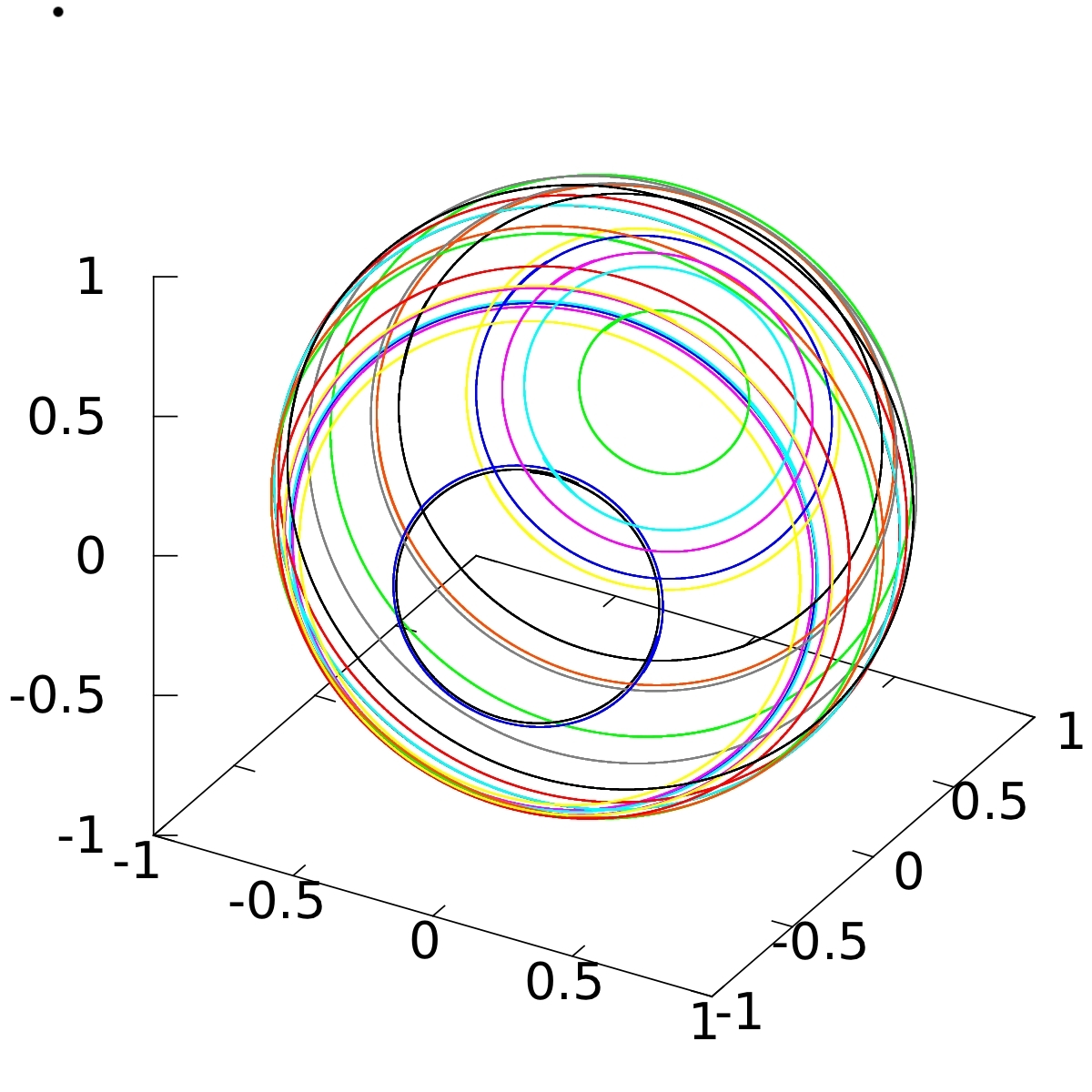}
    \caption{The trajectories of 25 agents with a randomly generated sign-symmetric interaction matrix starting from $t=0$ (left) and $t=0.03$ (right).}
    \label{fig:dancingeq}
\end{figure}

As mentioned in Section \ref{sec:app-soc}, the HK model was generalized in various ways. 
Since the model represents opinions, one may argue that more than one variable should be used to represent opinions on different issues. At the same time, opinions should belong to a compact set to be meaningful. Finally, a model with both features is obtained by projecting the HK model in dimension $d+1$ to the $d$-dimensional sphere
$\mathbb{S}^{d}$:
\begin{equation}\label{eq:HKsphere}
\dot x_{i} = \sum_{j=1}^{N} a_{ij} (x_{j} - \langle x_{i}, x_{j} \rangle x_{i}), \quad x_{i} \in \mathbb{S}^{d},
\end{equation}
where $a_{ij}$ are interaction weights: positive for attraction and negative for repulsion.
This model presents interesting dynamic features, including a new type of equilibria w.r.t. to the clusters of the HK model and new asymptotic states called
\emph{dancing equilibria}.
To describe the latter we introduce a piece of notation.
From \eqref{eq:HKsphere}, one can compute:
\begin{align}\label{eq:distances}
\frac{d}{dt}\langle x_{i},x_{j} \rangle =& \sum_{k\neq i} a_{ik} 
\left( \langle x_{k},x_{j} \rangle -\langle x_{k},x_{i} \rangle \langle x_{i},x_{j} \rangle \right)\nonumber\\
&+
\sum_{k\neq j} a_{jk} 
\left( \langle x_{k},x_{i} \rangle -\langle x_{k},x_{j} \rangle \langle x_{i},x_{j} \rangle \right),
\end{align}
where $\langle \cdot, \cdot \rangle$ is the scalar product in $\R^{d+1}$.
A trajectory along which the quantity \eqref{eq:distances} vanishes is called a \emph{dancing equilibrium}. 
The vanishing of 
\eqref{eq:distances} ensures that the relative position of agents on the sphere does not change. 
Therefore the trajectory could be obtained by fixing the agents at their initial position 
and rotating the whole sphere underneath. 
This behavior reminds us of cyclicity in many social systems and the visualization of such phenomenon on a two-dimensional sphere is provided in Figure \ref{fig:dancingeq}.
We refer the reader to \cite{CLP15} for details. Extensions to general manifolds are discussed in \cite{ACMPPRT17,AMP17}.

We notice that other interesting
asymptotic behaviors emerge in systems with many agents interacting in a social fashion. For instance, deep analysis was provided for specific configurations such as
consensus (or unique cluster), flocking (or alignment of velocities), and milling \cite{Degond2007,GomezSerrano2012,Motsch2014,MR2247927}.

\subsection{Time-evolving measures}
\label{sec:meas-evol}
\begin{figure}[h!]
    \centering
\includegraphics[height=5cm]{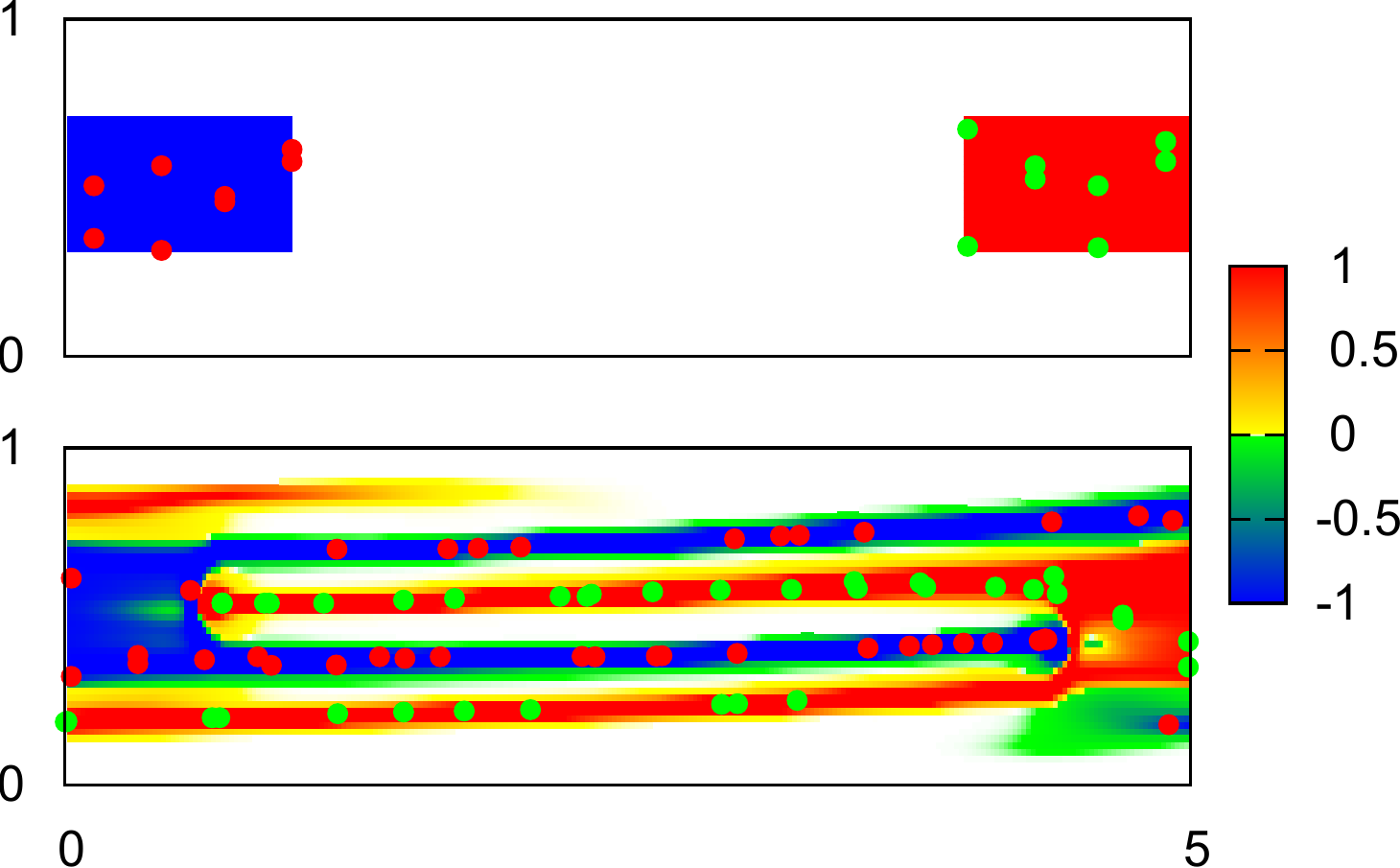}
    \caption{Simulation of a multiscale model for crossing flows.}
    \label{fig:cross-multi}
\end{figure}

Many of the patterns observed in pedestrian dynamics are significantly complex, such as, for example, the lane formation in opposite flows. The latter appears without initial self-organization and emerges as a form of symmetry breaking after the two flows start to interact. This happens ubiquitously at pedestrian crossings, see \cite{CPT14} for details.\\
The attempt to capture this complexity motivated the use of measure solutions instead of classical function solutions to dynamics modeling pedestrians, see
\cite{CPT11,PT09,PT11}.
In simple terms, a group of pedestrians occupying an environment $\Omega$ is represented by a measure $\mu$, which assigns to each (measurable) set $E\subset\Omega$ the number $\mu(E)$ quantifying pedestrians in $E$. This simple idea allows both continuum and discrete modeling, using either a $\mu$ absolutely continuous w.r.t. the Lebesgue measure (assuming $\Omega\subset\R^n$ for some $n$)
or an atomic measure.
Usually one assumed $\mu$ to be a Radon measure, but $\sigma$-additivity is the only necessary assumption to reflect the principle of additivity of the mass.
For $\mu$ to evolve in time,
one may use standard transport-type equations in the Eulerian setting such as:
\begin{equation}
	\frac{\partial\mu_t}{\partial t}+\nabla\cdot{(v\, \mu_t )}=0,
	\label{eq.cons.mass.strong}
\end{equation}
where $v [0,\,+\infty)\times\Omega\to\R^n$,
is the velocity at which the mass of pedestrians travels.
As usual for a measure, the equation has to be interpreted in the weak sense, i.e. 
for every test function $\phi$ and a.e. $t$ we require
\begin{equation}
	\frac{d}{dt}\int_{\R^d}\phi(x)\,d\mu_t(x)=\int_{\R^d}v(t,\,x)\cdot\nabla{\phi}(x)\,d\mu_t(x).
	\label{eq.cons.mass.weak}
\end{equation}
This model is of first-order, instead of second-order as the social force one, and it is convenient to represent the velocity $v$ as the sum of two terms:
\begin{equation}
		v(t,\,x)  =  v[\mu_t](x)
		 =  v_d(x)+v_i[\mu_t](x),
	\label{eq.v}
\end{equation}
where the dependence $v[\mu_t]$ is of functional type, i.e. $v$ is a map from the space of measures to the space of vector fields. The term $v_d$ represents a desired velocity, which depends only on the single pedestrian, and $v_i$ is an interaction term depending on the whole measure $\mu$.\\
This framework is well connected to mean-field approaches \cite{Golse16}, and it is convenient to use the Wasserstein distance \cite{villani}, which is defined as follows. Denote by $\mathcal{P}(\mathbb{R}^n)$ 
the space of probability measures, then for 
$\mu,\nu\in\mathbb{R}^n$ we set
\begin{equation}
\label{eq:W-def}
	W(\mu,\nu):= \inf\limits_{\pi\in \Pi(\mu,\nu)} \int\limits_{\R^n\times\R^n} |x-y| \mathrm{d}\pi(x,y),
\end{equation}
where $\Pi(\mu,\nu)$ indicates the set of transference plans between $\mu$ and $\nu$, i.e. the measures on 
$\R^n\times\R^n$ such that:
\begin{align*}
	\int\limits_{\R^n} \mathrm{d}\pi(\cdot, y)= \mathrm{d}\mu(\cdot),\quad \int\limits_{\R^n} \mathrm{d}\pi(x, \cdot)= \mathrm{d}\nu(\cdot).
\end{align*}  
The problem of finding the transference plan in \eqref{eq:W-def} is called the optimal transport problem. We refer the reader to \cite{villani,vi09} for a complete treatment of the subject and historical remarks.
In Figure \ref{fig:cross-multi}, we show a simulation of crossing flows with a multiscale measure model. The presence of both the discrete and continuous parts allows a symmetry breaking, thus generating the lanes observed in a pedestrian crossing. 

\section{New control problems and cost functions}
\label{sec:control-cost}





Besides new modeling approaches, the treatment of large-scale groups requires new ideas in the case of control problems.
In simple terms, the size of the group may render some classical approaches either unfeasible or not applicable by design.
For instance, in a real-world situation, one may be able to control only a small number of agents at the same time, or even at all.
An example of a new frontier in transportation is the problem of regulating traffic using autonomous vehicles from a given fleet, e.g. public buses. This corresponds to a control problem where the controllable agents are fixed and in small numbers. Moreover, communication issues may render it not possible to control in real time. 
In turn, mathematically, one has to consider a small number of fixed components of the system to be controlled. Most classical techniques may be not effective in such situations; for instance, a change of coordinates may not be feasible.
We illustrate some ideas to deal with such problems, including sparse controls, control across scales, and controls with bounded variation. 
We then illustrate new control problems, such as avoiding specific equilibria, and conclude with control problems for time-evolving measures.

\subsection{Sparse controls}
As mentioned above, the control of large systems requires the use of controls affecting a small number of agents. A useful concept is that of 
\emph{sparsity}.
In information theory \cite{Mallat09}, one looks for the most economical way to represent data,
which corresponds to the use of minimal
symbols from a dictionary, while,
in \emph{compressed sensing} \cite{FR13},
one looks for the representation of signals
using the minimal number of coefficients of a fixed base. 
Similarly, we define \emph{sparse controls} as the one using only a small number of agents. 

A convenient way to promote sparsity
is using $\ell^1$-type norms.
More precisely, consider $N$ agents, each evolving in $\R^d$, with an alignment dynamics as in \eqref{eq:CS}.
Define 
for every
$u,v \in \R^{Nd}$:
$$
B(u,v) = \frac{1}{2N^{2}}\sum_{i,j=1}^N \langle u_{i} - u_{j}, v_{i} - v_{j} \rangle = \frac{1}{N}  \sum_{i=1}^{N}\langle u_{i}, v_{i}\rangle - \langle \bar u, \bar v  \rangle,
$$
where $\langle \cdot, \cdot \rangle$ is the scalar product in $\R^{d}$.
The form $B$ is bilinear and can be used to define a variational problem to promote sparsity in the following way.
First fix a bound $M$ for the control action and $(x,v)\in \R^{Nd}\times \R^{Nd}$ 
(providing the positions and velocities of the agents). 
Then we define $U(x,v)$ as the set of solutions to: 
\begin{equation}\label{eq:variational}
\min \left( B (v, u) + \gamma(B(x,x)) \frac{1}{N} \sum_{i=1}^{N}\|u_{i}\| \right)  \quad \mbox{ subject to } \sum_{i=1}^{N}\|u_{i}\| \leq M\,,
\end{equation}
where
\begin{equation}\label{defgamma}
\gamma(b) = \int_{\sqrt{b}}^{\infty} a (\sqrt{2N}r) dr.
\end{equation}
The coefficient $\gamma$ measures the strength of future interaction and  if 
$B(v,v)<\gamma$ then the group tends to alignment without external intervention
\cite{HHK10}.
Thus the controls in $U(x,v)$ promote alignment, due to the term $B(u,v)$, and sparsity at the same time, due
to the $\ell^1$ terms in the cost and the constraint. Notice that $U(t,x)$ may be multivalued, thus generating a differential inclusion, but
sample-and-hold techniques can be used to select solutions with controls that are piecewise constant in time.
We refer the reader to \cite{CFPT13,CFPT15} for details.

\subsection{Controlling across scales and by leaders}\label{sec:across}
A classical tool in gas dynamics is passing to the mean-field limit in microscopic models \cite{Golse16}.
The limit consists of  partial differential equations of Vlasov or Boltzmann-type \cite{Perthame04}:
$$
\partial_t\mu+v \cdot \nabla_x \mu = \nabla_v \cdot \left [ \left ( H \star \mu \right ) \mu \right ]. 
$$
where $H$ is the interaction kernel, so the microscopic model reads 
$\dot{x}_i=v_i$, $\dot{v}_i = H\star\mu_N$,
with $\star$ the convolution and $\mu_N=\sum_i \delta_{(x_i,v_i)}$ the empirical measure. See \cite{CCH13} for applications of mean-field techniques to multi-agent systems.\\
In the controlled case, one would obtain the equation:
\begin{equation}\label{nufkup}
\partial_t\mu+v \cdot \nabla_x \mu = \nabla_v \cdot \left [ \left ( H \star \mu  \right ) \mu + \nu \right ], 
\end{equation}
but the limit of control policies $\nu$ may be a singular measure w.r.t. to $\mu$ thus having no effect on the dynamics.  Imaginatively, it is like trying to steer a river by means of toothpicks!\\
A possible solution is to restrict the control policies so that in the limit one obtains an equation of the type:
\begin{equation}\label{nufkup2}
\partial_t\mu+v \cdot \nabla_x \mu = \nabla_v \cdot \left [ \left ( H \star \mu  + f \right ) \mu \right ], 
\end{equation}
i.e. with $\nu$ absolutely continuous w.r.t. $\mu$ with Radon-Nikodym derivative
$\frac{d\nu}{d\mu}=f$. This approach was used in \cite{FS} assuming $f$ to be Lipschitz continuous in $(x,v)$.
While this is a severe restriction, it is possible to obtain a simultaneous $\Gamma$-limit and mean-field limit of
the optimal control problems.\\
An alternative idea is to split the population in two groups: leaders and followers, and assume that only the former can be controlled \cite{FPR14}. 
The complete microscopic dynamics reads:
\begin{equation}\label{eq:lead-follow-micro}
\begin{cases}
\dot y_k = w_k, & \\
\dot w_k =  H\star \mu_N(y_k,w_k) + H \star\mu_m(y_k,w_k)+u_k\, & k=1,\dots m,\\
\dot x_i = v_i, & \\
\dot v_i =  H\star \mu_N(x_i,v_i) + H \star\mu_m(x_i,v_i)\, & i=1,\dots N,
\end{cases}
\end{equation}
where $y_k$, respectively $w_k$, are the positions, respectively velocities, of the controlled leaders; $u_k$ are the controls; $x_i$, resp. $v_i$, are the positions, resp. velocities, of the followers;
and $\mu_m$, resp. $\mu_m$, is the empirical measure of followers, resp. leaders. We expect $m$ to be much smaller than $N$, so every control policy will be essentially sparse by design. In section \ref{sec:future}, we will consider specifically the control of traffic via autonomous vehicles. In that case, the ratio $\frac{m}{N}$ is called the penetration rate, and effective strategies should deal with the case of the penetration rate between 2\% and 5\%.\\
Optimal control problems for \eqref{eq:lead-follow-micro} can be state as follows:
\begin{equation}\label{eq:sparseoptcontr}
\min_{u \in L^1([0,T],U)} \int_{0}^T  \left \{ L(y(t),w(t),\mu_N(t))   +  \frac{1}{m} \sum_{k=1}^m |u_k(t)|  \right \}dt,
\end{equation}
where $U$ is the control set and $L$ the running cost or Lagrangian. The last term of $\ell^1$-type is further promoting sparsity for controls already acting only on leaders \cite{CCK12,CF65,WW11}.\\
In this setting, it is possible to define a mean-field limit of the controlled system letting $N\to\infty$ while keeping $m$ fixed. 
The resulting system is of mixed type with fully coupled controlled ODEs and a Vlasov-type PDE:
\begin{equation}
\begin{cases}
\dot y_k = w_k, & \\
\dot w_k =  H\star(\mu+\mu_m)(y_k,w_k)+u_k, & k=1,\dots , m, \quad t\in[0,T]\\
\partial_t\mu+v \cdot \nabla_x \mu = \nabla_v \cdot \left [ \left ( H \star (\mu +\mu_m) \right ) \mu \right ].
\end{cases}
\label{eq:lead-MF}
\end{equation}
Moreover, one can prove the  $\Gamma$-convergence
of the optimal control problem, so that
optimal controls of \eqref{eq:lead-follow-micro}-\eqref{eq:sparseoptcontr} converge to solutions to
the optimal control problem for \eqref{eq:lead-MF}
stated as:
\begin{equation}
\min_{u \in L^1([0,T],U)} \int_{0}^T  \left \{ L(y(t),w(t),\mu(t))   +  \frac{1}{m} \sum_{k=1}^m |u_k(t)|  \right \}dt.
\label{eq:MFcost}
\end{equation}
These results can be used in two different ways: designing or computing optimal controls at the microscopic scale \eqref{eq:lead-follow-micro}-\eqref{eq:sparseoptcontr} and proving properties of the limiting system, or using numerical methods for the mean-field scale \eqref{eq:lead-MF}-\eqref{eq:MFcost}
and project over the microscopic one for real applications. 

\subsection{Avoiding Black Swans}
In Section \ref{sec:sphere}, we showed how the behavior of large groups of agents may result in rich asymptotic regimes. At the same time, the classical stabilization problem to the origin may be replaced with stabilization problems to a periodic orbit or another limit set representing the asymptotic behavior of the group.
Often times such limit sets are locally attractive and Lyapunov-type techniques have been used, see \cite{CKRT22}. Here we present another innovative aspect: the avoidance of a limit set, called Black Swan. The latter, in the real socio-economic systems, represents a low-probability event with significant consequences for the system's dynamics, \cite{TalebBS}.

Following \cite{PPT19}, we focus on the general HK-type system:
\begin{equation}\label{eq:HK-gen}
\dot{x}_i=\frac{1}{N} \sum_j a(\|x_i-x_j\|) (x_j-x_i)+u_i(t),\quad
\sum_i \|u_i\|\leq M
\end{equation}
where $N$ is the number of agents and $u_i$ the controls.
To fix the ideas, assume that $a:]0,+\infty[\to ]0,+\infty[$ is locally Lipschitz with sublinear growth at infinity so that solutions to \eqref{eq:HK-gen} always exist (interpreting the case $x_i=x_j$ for some $i\not=j$ appropriately). 
The interest is to understand if there exist controls, possibly sparse, which may avoid consensus, i.e. the manifold ${\cal M}= \{x_1=\cdots=x_N\}$,
which represents the Black Swan. For instance, in a stock market consensus to sell may generate a financial collapse.
Four main results are stated in \cite{PPT19}:\\
{\bf Result 1.} If $\lim_{s\to 0} s\,a(s)=0$, i.e. $a$ doesn't blow up too fast for $s\to 0$, then for every $M$ and initial condition $x(0)$ there exist a (sparse) control strategy so that $x(t)\notin {\cal M}$ for every $t>0$. In simple words, using appropriate controls one can prevent the Black Swan to occur.\\
{\bf Result 2.} If $\lim_{s\to 0} s\,a(s)=+\infty$, i.e. $a$ blows up fast enough, then for every $M$ there exists a \emph{Black Hole}, i.e. a neighborhood ${\cal N}$ of ${\cal M}$ so that if $x(0)\in {\cal N}$, then the corresponding solution converges to ${\cal M}$ in finite time. In simple words, independently of the strength of the controls, if the system gets close enough to a consensus then the Black Swan is unavoidable.\\
{\bf Result 3.} If $\lim_{s\to +\infty} s\,a(s)=0$, i.e. $a$ decreases sufficiently fast to $0$ for $s\to 0$, then  for every $M$ there exists a safety region, i.e. a neighborhood ${\cal S}$ of infinity so that if $x(0)\in {\cal S}$, then there exists a (sparse) control keeping the evolution uniformly far away from the Black Hole.\\
{\bf Result 4.} Finally if $\lim_{s\to +\infty} s\,a(s)=+\infty$, i.e. $a$ does not decrease to $0$ fast enough for $s\to 0$, then  for every $M$ there exists basin of attraction, i.e. a neighborhood ${\cal B}$ of ${\cal M}$ which attracts trajectories from every initial data.\\
The results stated above were also generalized to the mean-field case, i.e. for the equation
\begin{equation}
\mu_t+\nabla\cdot \left( \left(\int a(\|x-y\|)(y-x)d\mu(x) + \chi_\omega u \right) \mu\right)=0
\end{equation}
where $\chi_\omega$ is the indicator function of the set $\omega$,
$u=u(x)$, $\|u\|_{L^\infty}\leq M$, and sparsity is modeled
as $\int_{\omega (t)}1\leq c$ for some fixed $c>0$.

\subsection{Bounded variation controls}
Another innovative approach to design \emph{parsimonious} controls is that of preventing chattering phenomena, i.e. the presence of an infinite number of switchings in optimal control problems. In \cite{Fuller}, the author considered the simple optimal control problem:
\begin{equation}\label{intro1}
\dot x_1=x_2,\ \dot x_2=u,\ u\in [-1,1],\quad 
\min \int_0^T x_1^2(t)\,dt,
\end{equation}
and showed that the unique optimal control satisfies:
$$
u(t)= 
\begin{cases}
1, & t\in(t_{2k},t_{2k+1}), \quad  k\in\N,\\
-1, & t\in(t_{2k+1},t_{2k+2}), \quad  k\in\N,
\end{cases}
$$
where $(t_k)$ is the increasing sequence of switching times, from $u=+1$ to $u=-1$ or viceversa. This phenomenon was shown to be generic for single-input systems in higher dimensions \cite{kupka}.

It is desirable to design control policies achieving costs close to the optimal one, but having a finite number of switchings. To do so, consider the general optimal control problem:
\begin{equation}\label{eq:OCP}
\begin{cases}
&\displaystyle{\min_{u\in{\cal U}}  \int_0^{t(u)}L(s,x(s),u(s))\,ds ,} \\
&\dot x=f(x,u),\quad u\in{\cal U},\\
&x(0)=x_0,\quad x(t(u))=0.
\end{cases}
\end{equation}
where ${\cal U}$ is the set of admissible controls steering the system to $0$ in time $t(u)$, $L$ the Lagrangian, and $f$ the dynamics. We assume that $f$ and $L$ are smooth, then
we approximate the optimal control problems with the family:
\begin{equation}\label{eq:OCP-eps}
\begin{cases}
&\displaystyle{\min_{u\in{\cal U}}  \left(\int_0^{t(u)}L(s,x(s),u(s))\,ds + \varepsilon\, \mathrm{TV}(u)\right),} \\
&\dot x=f(x,u),\quad u\in{\cal U},\\
&x(0)=x_0,\quad x(t(u))=0.
\end{cases}
\end{equation}
where $TV(u)$ indicates the total variation of $u$ (in time). Therefore, if a control has a finite cost, it must have a finite total variation so a finite number of switching. To state the first result,
we need some more notation. First, as usual
we indicate by $Lie_f$ the Lie algebra generated
by the vector field $\{f(\cdot,u): u\in U\}$
(where $U$ is the control set), and by 
$Lie_f(x)$ its evaluation at the point $x$.
We also recall that the system $\dot{x}=f(x,u)$ is
small-time locally controllable at $0$ if for any $\delta$
the set of points reachable from $0$ within time $\delta$ is a neighborhood of $0$. 
We are ready to state the following:
\begin{theorem}\label{th:th1-tv}
Consider the optimal control problem \eqref{eq:OCP}
and assume the following: i) for every $(t,x)$ the set $\{(f(x,u),L(t,x,u)+\gamma):u\in U, \gamma>0\}$ is convex; ii) $U$ is compact and $\max_{u(\cdot)\in{\cal U}} t(u)+\|x(\cdot,u(\cdot))\|<+\infty$, where $x(\cdot,u(\cdot))$ is the trajectory corresponding to $u(\cdot)$; iii) 
$Lie_f(0)={\cal F}=\R^n$ (where $x\in\R^n$); iv) and the system $\dot{x}=f(x,u)$ is small-time locally controllable at $0$. Then \eqref{eq:OCP-eps} admits solutions and: 
\begin{equation}\label{convcost}
\lim_{\varepsilon\to 0} \int_0^{t(u_\varepsilon)} L(t,x_\varepsilon(t),u_\varepsilon(t))\,dt=\int_0^{t(u^*)}L(t,x^*(t),u^*(t))\,dt ,
\end{equation}
where $x_\varepsilon(\cdot)$ is an optimal solution
to  \eqref{eq:OCP-eps} corresponding to the control $u_\varepsilon$. Moreover, $x_\varepsilon(\cdot)$ converges uniformly to an optimal solution of 
\eqref{eq:OCP}.
\end{theorem}
In simple words, Theorem \ref{th:th1-tv} states that every optimal control problem, which is sufficiently regular, admits approximations by control problems with a finite number of switchings. Under further regularity (see \cite{CGPT18} for details), one can achieve an estimate of the convergence rate:
\begin{equation*}
 \int_0^{t(u_\varepsilon)}L(t,x_{\varepsilon}(t),u_{\varepsilon}(t))\,dt-  \int_0^{t(u^*)} L(t,x^*(t),u^*(t))\,dt   =\mathrm{O}\left( \varepsilon^{\frac{\alpha\beta}{1+\alpha\beta}} \right),   
\end{equation*}
where $x^*$ is an optimal trajectory for \eqref{eq:OCP} corresponding to control $u^*$.
These techniques are useful also to approximate optimal controls for hybrid systems and for problems with state constraints, see \cite{CGPT18} for a detailed discussion and statements. 





\subsection{Control of time-evolving measures}
\label{sec:meas-control}
As mentioned in Section \ref{sec:meas-evol}, some problems are conveniently represented by time-evolving measures. The resulting transport-type equations have been studied in the framework of Wasserstein spaces, i.e. measure spaces endowed with Wassertein-type metrics.
For instance, consider the nonlinear transport equation:
\begin{equation}\label{eq:gen-tr}
\mu_t + \nabla\cdot (v[\mu]\mu)=0,
\end{equation}
where $v[\cdot ]$ highlight the functional dependence
of the velocity field $v$ on the measure $\mu$,
thus $v:{\cal P}\to {\cal C}(\R^n,\R^n)$ is a map
from the space of probability measures on $\R^n$, endowed with the Wasserstein distance $W$ \eqref{eq:W-def},  to the space of continuous vector fields on $\R^n$. In this setting, it is natural to state the conditions:
\begin{itemize}
    \item [(C1)] For every $\mu$, the vector field $v[\mu]$ is globally Lipschitz and bounded.
    \item[(C2)] The map $\mu\to v[\mu]$ is Lipschitz for the metric $W$ and the uniform norm
    $\|\cdot \|_\infty$.
\end{itemize}
These natural conditions allow the development of a fairly general framework allowing for a complete theory including approximation by Eulerian-type and Lagrangian-type numerical schemes, see
\cite{PR13}.\\
Extending this framework, one may consider control problems as follows. First consider a set-valued map $V: {\cal P}\rightrightarrows {\cal C}(\R^n,\R^n)$ associating to every probability measure a set of locally Lipschitz vector fields.
Then one can define the differential inclusions in Wassertein space:
\begin{equation}\label{eq:DI-W}
\mu_t\in -\nabla\cdot \left( V[\mu]\, \mu \right).
\end{equation}
As usual for differential inclusion problems, a solution to \eqref{eq:DI-W} is a curve $t\to \mu(t)$
for which there exists a measurable selection
$t\to w(t)\in V(\mu(t))$ so that \eqref{eq:gen-tr}
is satisfied with $v[\mu(t)]=w(t)$. We refer the reader to \cite{BF21} for the general theory and to \cite{BF21-1} for a version of Pontryagin Maximum Principle in this setting.\\
A different route to control problems for measures was taken in \cite{CMNP18,CMP17,CMP18,CMP18-1,CMP20},
where the authors considered an alternative formulation as follows. Starting from a classical controlled dynamics $\dot{x}=f(x,u)$, one can model a more general nonlocal controlled system as
\begin{equation}
\mu_t+\nabla\cdot (v(t)\mu)=0,
\end{equation}
where $v(t,x)\in \{f(x,u):u\in U\}$. This approach is more inherently connected to the microscopic dynamics $\dot{x}=f(x,u)$ than to the continuity equation
\eqref{eq:gen-tr}.\\
Going back to \eqref{eq:gen-tr}, it is quite straightforward to generalize the approach to measures with a total mass different from one (as for probability measures). However, some problems require the addition of source terms, which would modify the total mass in time, as in:
\begin{equation}
\mu_t+\nabla\cdot (v[\mu]\mu)=h[\mu].
\end{equation}
Then one has to modify the Wasserstein metric to compare measures with different masses. A possibility is to allow the addition/subtraction of mass at a fixed cost, which gives the metric:
\begin{equation}
W^{a,b}(\mu,\nu)=\inf_{\tilde{\mu}\leq\mu,\tilde{\nu}\leq \nu} \Big( a (|\mu-\tilde{mu}|+|\nu-\tilde{nu}|)
+b\,W(\mu,\nu)\Big)
\end{equation}
where $\tilde{\mu}\leq\mu$ means $\tilde{\mu}(A)\leq \mu(A)$ for every measurable set $A$ and
$|\cdot|$ is the total variation norm (or total mass).
In simple words, the distance $W^{a,b}$ considers all possible ways to arrange the mass of $\mu$ to the mass of $\nu$ by canceling mass with cost $a$ per unit of mass or transporting mass with cost $b$ per unit of mass.
This distance allows the development of a theory for equations with sources or sinks \cite{PR14,PR16}.
Alternative approaches based on optimal transport were proposed, and we refer the reader to \cite{PC19book} for a discussion including applications to image processing and artificial intelligence.

\section{Future perspectives}
\label{sec:future}

In this section, we provide some reflections on current and future perspectives. The new models and control problems and techniques described in Sections \ref{sec:new-model} and \ref{sec:control-cost}
posed new challenges. At the same time, they open new opportunities for impacting different areas of engineering and social sciences.
First, we present a current view on the limiting approaches for dealing with large systems. This includes a new type of equation for measures as well as a new type of limit.
Then we focus on a research effort in vehicular traffic, which gave rise to a large-scale real-world application. More precisely, we illustrate the work of a consortium, called CIRCLES \cite{circles}, on smoothing traffic with autonomous vehicles (briefly AVs). Starting in 2014 a group of researchers showed how a small number of AVs can impact traffic by dampening stop-and-go waves, reaching the peak of a 100 AV experiment performed in November 2022 on the I-24 (a multilane interstate highway) in Nashville, TN.

\subsection{A general view  on equations for measures and limits}
Building upon the modeling effort illustrated in Section \ref{sec:meas-control}, recently a new type of equation was proposed, called Measure Differential Equations (briefly MDE), see \cite{piccoli2019}. The main idea is to generalize the concept of ordinary differential equations on a manifold to measures. 
Let us start by noticing that an ODE is defined by assigning a vector field $V:M\to TM$,
$x\to (x,v(x))$, assigning to every point $x\in M$
a tangent vector of the tangent bundle $TM$.
Specifically, the map definition $x\to (x,v(x))$
forces the tangent vector associated to $x$ to belong to the tangent space $T_xM$ of $M$ at the point $x$, i.e. to the fiber of $TM$ over $x$.
An alternative way to state this is to require that $\pi(v(x))=x$, where $\pi:TM\to M$ is the projection over the base, i.e. $\pi(x,v)=x$.
A the natural immersion of $M$ into ${\cal P}(M)$,
the space of probability measures over $M$, is given by $x\to\delta_x$ associating to every point the Dirac delta centered at the point. 
This immersion can be defined for every manifold,
thus can be used to extend the concept of vector field to that of Probability Vector Field (briefly PVF) on $M$, which is a map $V:{\cal P}(M)\to {\cal P}(TM)$ associating to every probability measure on $M$, a probability measure on $TM$, with the only requirement
$\pi\# V(\mu)=\mu$, where $\#$ is the push-forward of measures. Given a PVF $V$, we can define an MDE as:
\begin{equation}\label{eq:MDE}
\dot{\mu}=V[\mu]
\end{equation}
where a solution is interpreted in weak sense. This means that for every test function $f$ we have:
\begin{equation}\label{eq:MDE-sol}
\frac{d}{dt}\int_{x\in M} f\,d[\mu(t)](x) =
\int_{TM}(\nabla_x f\cdot v) \, dV[\mu](t,x).
\end{equation}
Existence of solutions can be proved using standard assumptions on $V$ of sublinear growth (of measure support) and continuity w.r.t. to Wasserstein metrics on $M$ and $TM$. However, uniqueness can be recovered only at the level of semigroup using Lipschitz-type assumptions based on the functional:
\begin{equation}
{\cal W}(V,W)=\min\{
\int_{TM}|v-w|dT(x,v,y,w): T\in\Pi (V,W),
\pi_{13}\#T \ optimal\}
\end{equation}
where $V,W\in {\cal P}(TM)$, $\Pi$ is the set of transference plans and $\pi_{13}:TM\times TM\to M\times M$ is the projection over first and third component. In simple words, we consider all transference plans between measures on $TM$, with the requirement that the projection on the bases to be optimal (in the sense of optimal transport). 
This approach was recently generalized considering
the one-sided Lipschitz condition for ${\cal W}$
and making ingenious connections with evolution variational inequalities \cite{CSS22}. We notice that this approach allows for uniqueness at the level of solutions, thus opening the door to control problems. On the other side, if one wants to use a non-variational interpretation, then the setting of MDEs suffers from the fact that PVFs enjoy only a structure of abelian monoid and not of Lie algebra. A possibility is to explore control properties at the level of the flows, instead of solutions, taking advantage of the geometric control techniques \cite{ASbook}.

\begin{figure}[h]
\begin{center}
\resizebox{12cm}{!}{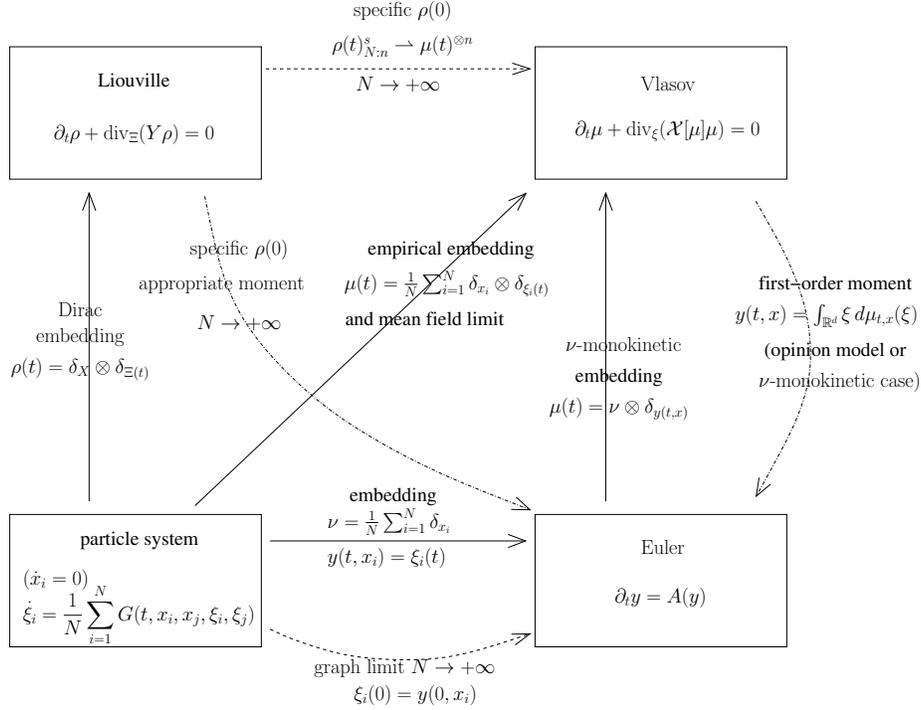}
\end{center}
\caption{Relationships between particle (microscopic) system, Liouville (probabilistic) equation, Vlasov (mesoscopic, mean-field) equation, Euler (macroscopic, graph limit) equation. From \cite{PT22} with permission from the authors.}\label{fig_embeddings}
\end{figure}

Another intriguing vision was given by the recent preprint \cite{PT22}, where multiple scales and concepts of limits were discussed. With permission from the authors, we report a figure representing the general picture of limiting processes, see Figure \ref{fig_embeddings}.
The authors start with a general microscopic system of the type:
\begin{equation}\label{eq:gen-sys-PT}
\dot{\xi}_i=\frac{1}{N} \sum_j G(t,x_i,x_j,\xi_i,\xi_j)
\end{equation}
where $x_i$ represent parameters describing the $i$-th agent, while $\xi$ is its state.
Comparing with systems presented in Section \ref{sec:app}, the main difference is the presence of the variable $x_i$, which is assumed to be constant in time, thus $\dot{x}_i=0$. This allows the interaction function $G$ to depend on the specific particle bypassing the classical
limitation of \emph{indistinguishibility} of particles assumed in mean-field theories \cite{Golse16}.\\
An alternative way to send the number of agents $N$ to infinity is via the so-called \emph{graph-limit}
by interpreting the right-hand side of \eqref{eq:gen-sys-PT} as a Riemann sum. Another interpretation is to look at the graph of interacting agents as a \emph{graphon}, i.e. a function from $[0,1]^2$ to $[0,1]$, and passing to the limit in the cut norm \cite{Med14}. For recent results see also \cite{AD21,BCZ19}.\\
In \cite{PT22} the authors show other connections via mean-field, graph-limit, and embeddings. A future direction of great interest will be to understand how control problems behave with respect to such limits, as in the spirit of Section \ref{sec:across}.\\
Along the same notes, we want to point to another recent contribution
\cite{CLOS22}, where the authors, besides the classical Lagrangian and Eulerian point of view, consider the \emph{Kantorovich} formulation. The latter is based on representing solutions to equations as measures over the space of continuous curves. Using the superposition principle \cite{AGS08}, one is able to connect the Lagrangian and Eulerian formulation to the Kantorovich one. In this framework, the authors present various equivalence and $\Gamma$-convergence results, but it would be interesting to understand sparsity, parsimonious controls, and necessary and sufficient conditions for optimality.

\subsection{Real world applications}
In this section, we report experimental results concerning the control of traffic via autonomous vehicles (briefly AVs). More precisely, we first discuss the problem of stop-and-go waves and their impact on traffic. Then, we show how a single AV can dissipate waves on a ring road with other 20 cars. Finally, we describe a recent traffic experiment with 100 AVs on an open highway: the I-24 in Nashville, Tennessee. The latter is, at the present time, the largest traffic experiment ever performed with AVs on an open highway.

\subsubsection{Stop-and-go waves and fuel consumption}
As recalled in Section \ref{sec:app-traffic}, stop-and-go waves may appear ubiquitously due to differences in human driving, also in a confined setting
\cite{ sugiyama2008,treiterer1974}. They can be explained in terms of instabilities of traffic flow \cite{cui2017stabilizing}, and their impact includes increasing fuel consumption and increased breaking \cite{RamadanSeibold2017}. Several studies proposed strategies
to smooth traffic via AVs in simulation
\cite{PhysRevE.69.066110,GUERIAU2016266,TALEBPOUR2016143,7362183}. 

\begin{figure}
    \centering
    \includegraphics[width=\textwidth]{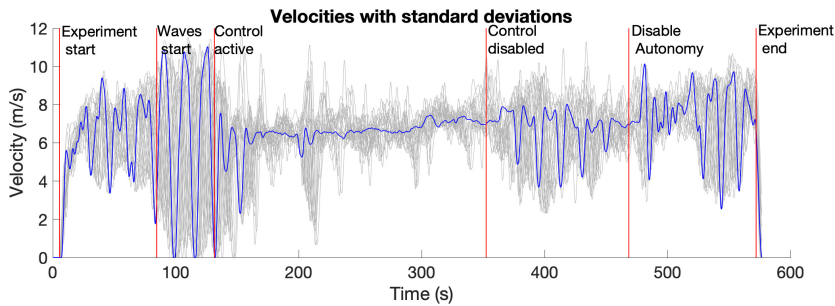}
    \caption{Trajectories of the field experiment on a ring-road: AV trajectory in blue, human-driven vehicles' trajectories in gray.}
    \label{fig:ch5:exp_velocities}
\end{figure}

A first set of experiments to dissipate waves via a single AV were performed on a ring-road \cite{STERN2019351,stern2018dissipation,WU201982},
which is the same setting as the seminal experiment \cite{ sugiyama2008}. 
A single AV with level 2 autonomy (according to the SAE taxonomy) out of 22 cars was able to dissipate waves, improve fuel economy, and reduce significant breaking events. 
The experiment started with the same setting as the Japanese one,
with the driver instructed only to keep a safe distance and the AV being human-driven.
Figure \ref{fig:ch5:exp_velocities} reports the velocities of all human-driven vehicles (in gray) and the AV one (in blue).
One can notice the appearance of a wave, which, after around 90 seconds,
becomes stable causing all vehicles to accelerate and brake periodically.
The velocity oscillations are significant ranging between $0$ and $11$ m/s. Then the AV automatic control is activated (the human commands only the steering wheel) and the wave is dissipated in a short time.  
At the end of the experiment the AV turns back to human-driven mode and the wave reappears. The final result is a reduction of fuel consumption
of $43\%$ (measured with OBD-II on all vehicles) and a reduction of heavy breaking events (decrease of more than 1m/s in 1s) by $98\%$.

\subsubsection{The I-24 MOTION testbed} 
\begin{figure}
    \centering
    \includegraphics[width = 0.465\textwidth]{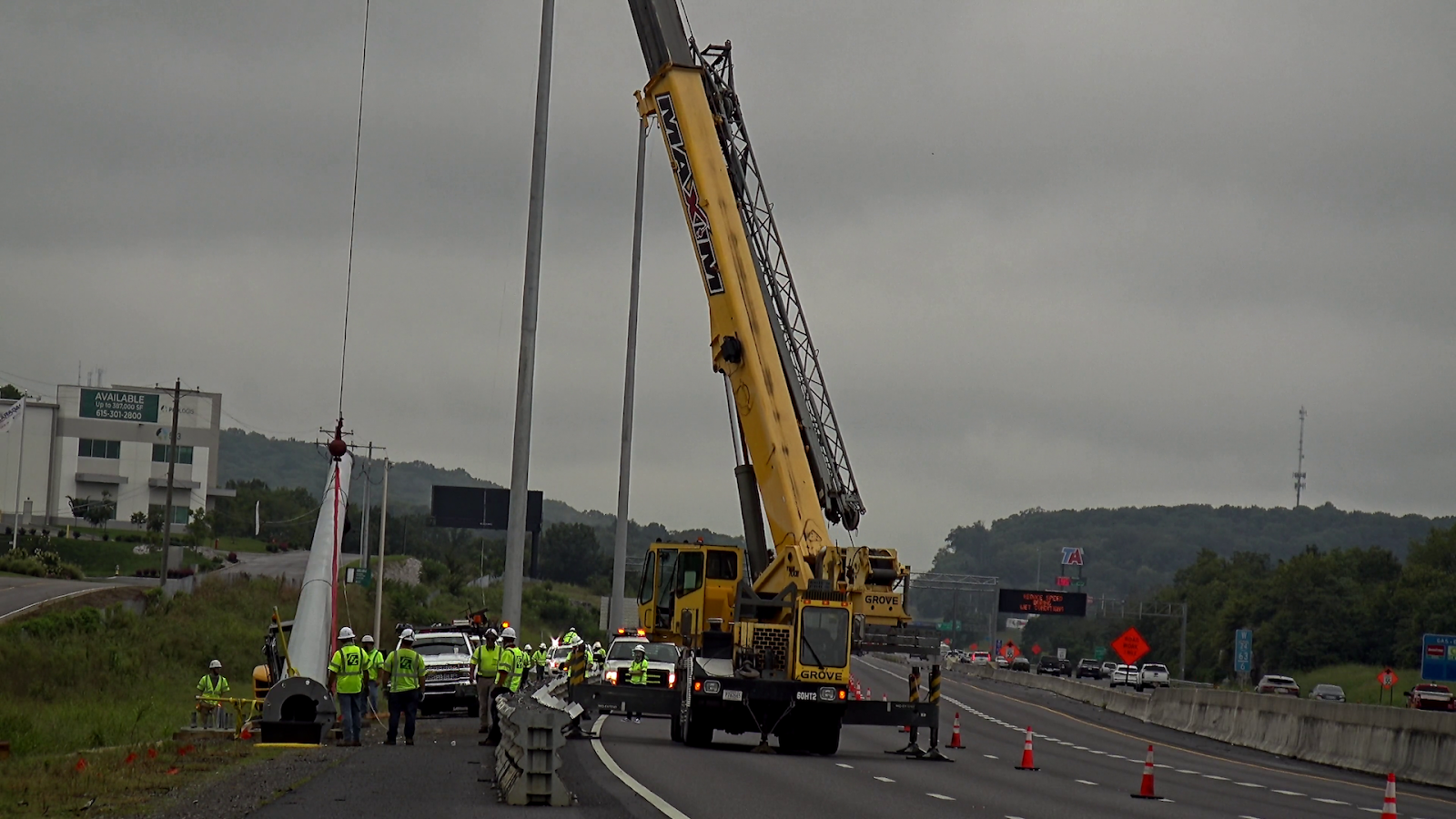}
    \includegraphics[width = 0.52\textwidth]{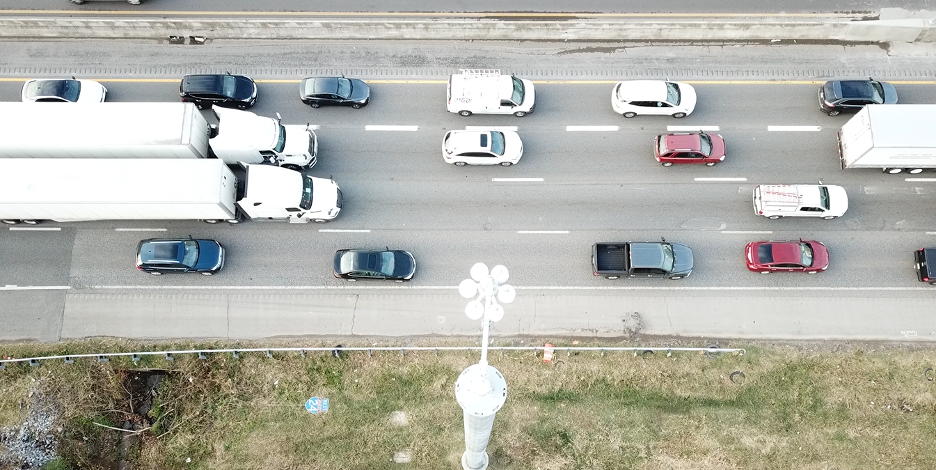}
    \caption{Left: Installation of a roadside pole. Right: aerial view of the installed pole with 8 cameras  and I-24 traffic underneath.}
    \label{fig:roadsidepoles}
\end{figure}

A sizeable testbed was developed on the interstate I-24 in Nashville, TN,
called I-24 MOTION\cite{I-24website}. This testbed was completed in October 2022  by the Tennessee Department of Transportation and Vanderbilt University with almost 300 4K cameras mounted on more than forty \(110\)ft roadside poles.
Each pole hosts six cameras to allow complete coverage of 500 linear feet of highway. The system captures more than 150,000 vehicles daily \cite{gloudemans2020interstate}. For comparison, each day the system captures data that are three orders of magnitude larger than the celebrated NGSIM data set \cite{NGSIMurl}.
Figure \ref{fig:roadsidepoles} show the construction site for the installation of one of the poles. 

\subsubsection{The 100-AV experiment}
\begin{figure}[h!]
    \centering
    \includegraphics[trim={1cm 0cm 2.6cm 0}, clip,scale=0.33]{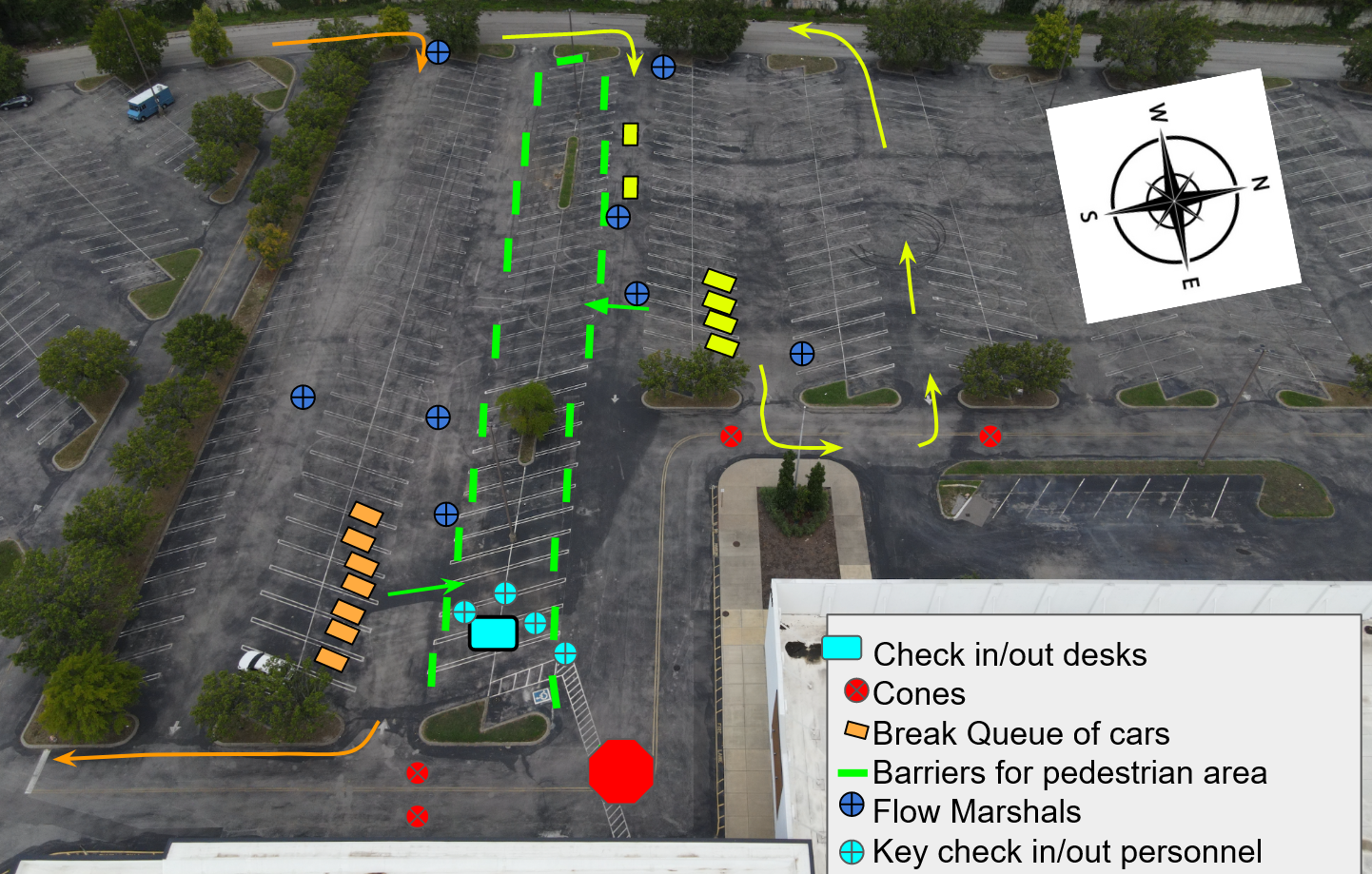}
    \includegraphics[trim={3cm 6cm 1cm 4cm 0}, clip, scale=0.81]{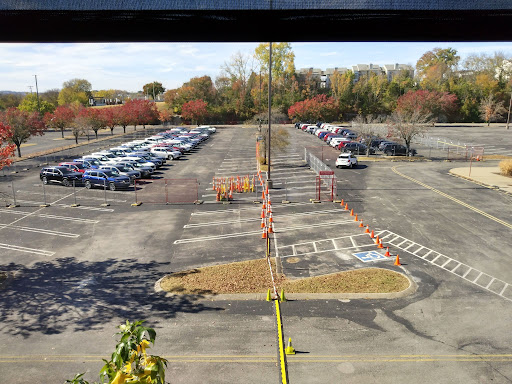}
    
    \caption{Top: A plan of the parking lot from which the cars left to reach the I-24. Orange and yellow routes are highlighted. Bottom: A photo of the lot during the experiment with the two groups of cars (corresponding to the orange and yellow routes).}
    \label{fig:Lot layout}
\end{figure}

In November 2022, the CIRCLES consortium \cite{circles} performed a 100 AV experiment on the I-24 testbed. The research team was comprised of more than 60 researchers from the CIRCLES main partners (UC Berkeley, Rutgers-Camden, Temple, Vanderbilt) as well as others (Arizona, ENS, and others).
Moreover, three OEMs participated: Nissan, GM, and Toyota. 
During the week of November 14th, 2022, the team deployed the vehicles, again with level 2 autonomy, during the morning rush hour for five days testing several controls developed with mathematical model techniques as well as AI tools. More than 150 drivers were trained and sent to the highways on two fixed routes around the portion of the highway monitored by cameras. The experiment was planned for several months to allow the safe deployment of vehicles and an appropriate penetration (percentage of controlled vehicles on overall traffic). Figure\ref{fig:Lot layout} depicts the planned operations on the parking lot, from which cars left to reach the I-24, and shows a picture from one day of the experiment.

\subsection{Future experiments}
The 100 AVs experiment represented the largest experiment worldwide of its kind. As we are writing, the massive experimental data are analyzed to quantify the impact of AV regulation on traffic. Preliminary evidence showed a significant impact.\\
The experiment showed evidence of the feasibility of interventions at scale on real traffic. Future development may include the following:
\begin{itemize}
    \item Traffic regulation via fleets, e.g. buses, taxis, and others.
    \item Recruitment of community stakeholders, such as commuters, to reduce congestion and fuel consumption.
    \item Exportation of the ideas to other areas, including pedestrian crowds, robotic formations, and socio-economic systems. 
\end{itemize}

\section*{Acknowledgements}
This material is based upon work supported by the National Science Foundation under the grants CNS-1837481 and CNS-1446715. The work of the author was supported by the U.S. Department of Energy’s Office of Energy Efficiency and Renewable Energy (EERE) under the Vehicle Technologies Office award number CID DE-EE0008872. The views expressed herein do not necessarily represent the views of the U.S. Department of Energy or the United States Government.
The author thanks Christopher Denaro, Sean McQuade, Ryan Weightman, for discussions and proofreading of the paper.

\small
\bibliographystyle{abbrv}
\bibliography{CrowdsReview,Full,Biblio}

\begin{thebibliography}{100}

\bibitem{I-24website}
I-24 {MOTION}.
\newblock \url{https://i24motion.org/}.

\bibitem{NGSIMurl}
{N}ext {G}eneration {SIM}ulation.
\newblock http://ngsim-community.org/, 2013.

\bibitem{ASbook}
A.~A. Agrachev and Y.~L. Sachkov.
\newblock {\em Control theory from the geometric viewpoint}, volume~87 of {\em
  Encyclopaedia of Mathematical Sciences}.
\newblock Springer-Verlag, Berlin, 2004.
\newblock Control Theory and Optimization, II.

\bibitem{IDM22}
S.~Albeaik, A.~Bayen, M.~T. Chiri, X.~Gong, A.~Hayat, N.~Kardous, A.~Keimer,
  S.~T. McQuade, B.~Piccoli, and Y.~You.
\newblock Limitations and improvements of the intelligent driver model ({IDM}).
\newblock {\em SIAM J. Appl. Dyn. Syst.}, 21(3):1862--1892, 2022.

\bibitem{AGS08}
L.~Ambrosio, N.~Gigli, and G.~Savar\'{e}.
\newblock {\em Gradient flows in metric spaces and in the space of probability
  measures}.
\newblock Lectures in Mathematics ETH Z\"{u}rich. Birkh\"{a}user Verlag, Basel,
  second edition, 2008.

\bibitem{Arechavaleta2008}
G.~Arechavaleta, J.-P. Laumond, H.~Hicheur, and A.~Berthoz.
\newblock An optimality principle governing human walking.
\newblock {\em IEEE Transactions on Robotics}, 24:5--14, 2008.

\bibitem{AwRascle}
A.~Aw and M.~Rascle.
\newblock Resurection of second order models of traffic flow.
\newblock {\em SIAM J. Appl. Math.}, 60:916--944, 2000.

\bibitem{ACMPPRT17}
A.~Aydo\u{g}du, M.~Caponigro, S.~McQuade, B.~Piccoli, N.~Pouradier~Duteil,
  F.~Rossi, and E.~Tr\'{e}lat.
\newblock Interaction network, state space, and control in social dynamics.
\newblock In {\em Active particles. {V}ol. 1. {A}dvances in theory, models, and
  applications}, Model. Simul. Sci. Eng. Technol., pages 99--140.
  Birkh\"{a}user/Springer, Cham, 2017.

\bibitem{AMP17}
A.~Aydo\u{g}du, S.~T. McQuade, and N.~Pouradier~Duteil.
\newblock Opinion dynamics on a general compact {R}iemannian manifold.
\newblock {\em Netw. Heterog. Media}, 12(3):489--523, 2017.

\bibitem{AD21}
N.~Ayi and N.~{Pouradier Duteil}.
\newblock Mean-field and graph limits for collective dynamics models with
  time-varying weights.
\newblock {\em Journal of Differential Equations}, 299:65--110, 2021.

\bibitem{bellomo20review}
B.~Aylaj, N.~Bellomo, L.~Gibelli, and D.~Knopoff.
\newblock Crowd dynamics by kinetic theory modeling: Complexity, modeling,
  simulations, and safety, 2020.

\bibitem{B08}
M.~Ballerini, N.~Cabibbo, R.~Candelier, A.~Cavagna, E.~Cisbani, I.~Giardina,
  V.~Lecomte, A.~Orlandi, G.~Parisi, A.~Procaccini, M.~Viale, and
  V.~Zdravkovic.
\newblock Interaction ruling animal collective behavior depends on topological
  rather than metric distance: Evidence from a field study.
\newblock {\em Proceedings of the National Academy of Sciences},
  105(4):1232--1237, 2008.

\bibitem{Bando1995}
M.~Bando, K.~Hasebe, A.~Nakayama, A.~Shibata, and Y.~Sugiyama.
\newblock Dynamical model of traffic congestion and numerical simulation.
\newblock {\em Physical review E}, 51(2):1035, 1995.

\bibitem{BDMGPbook}
A.~Bayen, M.~Delle~Monache, M.~Garavello, P.~Goatin, and B.~Piccoli.
\newblock {\em Control Problems for Conservation Laws with Traffic
  Applications}.
\newblock Springer Nature (Birkhauser), Cham, 2022.

\bibitem{circles}
A.~Bayen, J.~Lee, B.~Seibold, B.~Piccoli, J.~Sprinkle, and D.~Work.
\newblock Circles consortium project website.
\newblock \url{https://circles-consortium.github.io}, 2021.

\bibitem{Bellomo2016}
N.~Bellomo, D.~Clarke, L.~Gibelli, P.~Townsend, and B.~Vreugdenhil.
\newblock Human behaviours in evacuation crowd dynamics: From modelling to
  “big data” toward crisis management.
\newblock {\em Physics of Life Reviews}, 18:1 -- 21, 2016.

\bibitem{MR2974186}
N.~Bellomo and J.~Soler.
\newblock On the mathematical theory of the dynamics of swarms viewed as
  complex systems.
\newblock {\em Math. Models Methods Appl. Sci.}, 22(suppl. 1):1140006, 29,
  2012.

\bibitem{BCZ19}
U.~Biccari, D.~Ko, and E.~Zuazua.
\newblock Dynamics and control for multi-agent networked systems: A
  finite-difference approach.
\newblock {\em Mathematical Models and Methods in Applied Sciences},
  29(04):755--790, 2019.

\bibitem{BF21}
B.~Bonnet and H.~Frankowska.
\newblock Differential inclusions in {W}asserstein spaces: the
  {C}auchy-{L}ipschitz framework.
\newblock {\em J. Differential Equations}, 271:594--637, 2021.

\bibitem{BF21-1}
B.~Bonnet and H.~Frankowska.
\newblock Necessary optimality conditions for optimal control problems in
  {W}asserstein spaces.
\newblock {\em Appl. Math. Optim.}, 84(suppl. 2):S1281--S1330, 2021.

\bibitem{BP07}
A.~Bressan and B.~Piccoli.
\newblock {\em Introduction to the mathematical theory of control}, volume~2 of
  {\em AIMS Series on Applied Mathematics}.
\newblock AIMS, Springfield, MO, 2007.

\bibitem{CFPT13}
M.~Caponigro, M.~Fornasier, B.~Piccoli, and E.~Tr\'{e}lat.
\newblock Sparse stabilization and optimal control of the {C}ucker-{S}male
  model.
\newblock {\em Math. Control Relat. Fields}, 3(4):447--466, 2013.

\bibitem{CFPT15}
M.~Caponigro, M.~Fornasier, B.~Piccoli, and E.~Tr{\'e}lat.
\newblock Sparse stabilization and control of alignment models.
\newblock {\em Mathematical Models and Methods in Applied Sciences},
  25(3):521--564, 2015.

\bibitem{CGPT18}
M.~Caponigro, R.~Ghezzi, B.~Piccoli, and E.~Tr\'{e}lat.
\newblock Regularization of chattering phenomena via bounded variation
  controls.
\newblock {\em IEEE Trans. Automat. Control}, 63(7):2046--2060, 2018.

\bibitem{CLP15}
M.~Caponigro, A.~C. Lai, and B.~Piccoli.
\newblock A nonlinear model of opinion formation on the sphere.
\newblock {\em Discrete and Continuous Dynamical Systems Ser. A},
  (9):4241--4268, 2015.

\bibitem{CCH13}
J.~A. Carrillo, Y.-P. Choi, and M.~Hauray.
\newblock {The derivation of swarming models: mean-field limit and Wasserstein
  distances}.
\newblock In A.~Muntean and F.~Toschi, editors, {\em Collective Dynamics from
  Bacteria to Crowds}, CISM International Centre for Mechanical Sciences, pages
  1--46. Springer, 2014.

\bibitem{CKRT22}
J.~A. Carrillo, D.~Kalise, F.~Rossi, and E.~Tr\'{e}lat.
\newblock Controlling swarms toward flocks and mills.
\newblock {\em SIAM J. Control Optim.}, 60(3):1863--1891, 2022.

\bibitem{CCK12}
E.~Casas, C.~Clason, and K.~Kunisch.
\newblock Approximation of elliptic control problems in measure spaces with
  sparse solutions.
\newblock {\em SIAM J. Control Optim.}, 50(4):1735--1752, 2012.

\bibitem{CLOS22}
G.~Cavagnari, S.~Lisini, C.~Orrieri, and G.~Savaré.
\newblock Lagrangian, eulerian and kantorovich formulations of multi-agent
  optimal control problems: Equivalence and gamma-convergence.
\newblock {\em Journal of Differential Equations}, 322:268--364, 2022.

\bibitem{CMNP18}
G.~Cavagnari, A.~Marigonda, K.~T. Nguyen, and F.~S. Priuli.
\newblock Generalized control systems in the space of probability measures.
\newblock {\em Set-Valued and Variational Analysis}, 26(3):663--691, 2018.

\bibitem{CMP17}
G.~Cavagnari, A.~Marigonda, and B.~Piccoli.
\newblock Optimal synchronization problem for a multi-agent system.
\newblock {\em Netw. Heterog. Media}, 12(2):277--295, 2017.

\bibitem{CMP18-1}
G.~Cavagnari, A.~Marigonda, and B.~Piccoli.
\newblock Averaged time-optimal control problem in the space of positive
  {B}orel measures.
\newblock {\em ESAIM Control Optim. Calc. Var.}, 24(2):721--740, 2018.

\bibitem{CMP18}
G.~Cavagnari, A.~Marigonda, and B.~Piccoli.
\newblock Superposition principle for differential inclusions.
\newblock In I.~Lirkov and S.~Margenov, editors, {\em Large-Scale Scientific
  Computing}, pages 201--209, Cham, 2018. Springer International Publishing.

\bibitem{CMP20}
G.~Cavagnari, A.~Marigonda, and B.~Piccoli.
\newblock Generalized dynamic programming principle and sparse mean-field
  control problems.
\newblock {\em Journal of Mathematical Analysis and Applications},
  481(1):123437, 2020.

\bibitem{CSS22}
G.~Cavagnari, G.~Savar{\'e}, and G.~E. Sodini.
\newblock Dissipative probability vector fields and generation of evolution
  semigroups in wasserstein spaces.
\newblock {\em Probability Theory and Related Fields}, 2022.

\bibitem{CFPR21}
F.~Ceragioli, P.~Frasca, B.~Piccoli, and F.~Rossi.
\newblock Generalized solutions to opinion dynamics models with
  discontinuities, 2021.

\bibitem{Cercignani1994}
C.~Cercignani, R.~Illner, and M.~Pulvirenti.
\newblock {\em The mathematical theory of dilute gases}, volume 106 of {\em
  Applied Mathematical Sciences}.
\newblock Springer-Verlag, New York, 1994.

\bibitem{Chalons2017}
C.~Chalons, M.~L. Delle~Monache, and P.~Goatin.
\newblock A conservative scheme for non-classical solutions to a strongly
  coupled {PDE}-{ODE} problem.
\newblock {\em Interfaces Free Bound.}, 19(4):553--570, 2017.

\bibitem{Chitour2012}
Y.~Chitour, F.~Jean, and P.~Mason.
\newblock Optimal control models of goal-oriented human locomotion.
\newblock {\em SIAM Journal on Control and Optimization}, 50(1):147--170, 2012.

\bibitem{MR2165531}
S.~Cordier, L.~Pareschi, and G.~Toscani.
\newblock On a kinetic model for a simple market economy.
\newblock {\em J. Stat. Phys.}, 120(1-2):253--277, 2005.

\bibitem{CF65}
A.~J. Craig and I.~Flugge-Lotz.
\newblock {Investigation of Optimal Control With a Minimum-Fuel Consumption
  Criterion for a Fourth-Order Plant With Two Control Inputs; Synthesis of an
  Efficient Suboptimal Control}.
\newblock {\em Journal of Basic Engineering}, 87(1):39--57, 03 1965.

\bibitem{CFP11}
E.~Cristiani, P.~Frasca, and B.~Piccoli.
\newblock Effects of anisotropic interactions on the structure of animal
  groups.
\newblock {\em J. Math. Biol.}, 62(4):569--588, 2011.

\bibitem{CPTbook}
E.~Cristiani, B.~Piccoli, and A.~Tosin.
\newblock Multiscale modeling of granular flows with application to crowd
  dynamics.
\newblock {\em Multiscale Model. Simul.}, 9(1):155--182, 2011.

\bibitem{CPT11}
E.~Cristiani, B.~Piccoli, and A.~Tosin.
\newblock Multiscale modeling of granular flows with application to crowd
  dynamics.
\newblock {\em Multiscale Model. Simul.}, 9(1):155--182, 2011.

\bibitem{Cristiani2014}
E.~Cristiani, B.~Piccoli, and A.~Tosin.
\newblock {\em Multiscale Modeling of Pedestrian Dynamics}, volume~12 of {\em
  MS\&A. Model. Simul. Appl.}
\newblock Springer MS \& A: Modeling, Simulation and Applications, 2014.

\bibitem{CPT14}
E.~Cristiani, B.~Piccoli, and A.~Tosin.
\newblock {\em Multiscale Modeling of Pedestrian Dynamics}, volume~12.
\newblock Springer MS \& A: Modeling, Simulation and Applications, 2014.

\bibitem{Cucker2007}
F.~Cucker and S.~Smale.
\newblock Emergent behavior in flocks.
\newblock {\em IEEE Transactions on Automatic Control}, 52:852--862, 2007.

\bibitem{cui2017stabilizing}
S.~Cui, B.~Seibold, R.~Stern, and D.~B. Work.
\newblock Stabilizing traffic flow via a single autonomous vehicle:
  Possibilities and limitations.
\newblock In {\em 2017 IEEE Intelligent Vehicles Symposium (IV)}, pages
  1336--1341. IEEE, 2017.

\bibitem{DaganzoCritic}
C.~F. Daganzo.
\newblock Requiem for second-order fluid approximation to traffic flow.
\newblock {\em Transportation Research Part B}, 29(4):277--286, 1995.

\bibitem{DalMasoBook}
G.~Dal~Maso.
\newblock {\em An introduction to $\Gamma$ convergence}.
\newblock Springer, New York, 1993.

\bibitem{PhysRevE.69.066110}
L.~C. Davis.
\newblock Effect of adaptive cruise control systems on traffic flow.
\newblock {\em Phys. Rev. E}, 69:066110, Jun 2004.

\bibitem{Degond2013a}
P.~Degond, J.-G. Liu, S.~Motsch, and V.~Panferov.
\newblock Hydrodynamic models of self-organized dynamics: derivation and
  existence theory.
\newblock {\em Methods Appl. Anal.}, 20(2):89--114, 2013.

\bibitem{Degond2007}
P.~Degond and S.~Motsch.
\newblock Macroscopic limit of self-driven particles with orientation
  interaction.
\newblock {\em C. R. Math. Acad. Sci. Paris}, 345(10):555--560, 2007.

\bibitem{DMG2014num}
M.~L. Delle~Monache and P.~Goatin.
\newblock A front tracking method for a strongly coupled {PDE}-{ODE} system
  with moving density constraints in traffic flow.
\newblock {\em Discrete Contin. Dyn. Syst. Ser. S}, 7(3):435--447, 2014.

\bibitem{DMG2014exist}
M.~L. Delle~Monache and P.~Goatin.
\newblock Scalar conservation laws with moving constraints arising in traffic
  flow modeling: an existence result.
\newblock {\em J. Differential Equations}, 257(11):4015--4029, 2014.

\bibitem{Farina2016}
F.~Farina, D.~Fontanelli, A.~Garulli, A.~Giannitrapani, and D.~Prattichizzo.
\newblock When {H}elbing meets {L}aumond: {T}he {H}eaded {S}ocial {F}orce
  {M}odel.
\newblock In {\em 2016 IEEE 55th Conference on Decision and Control (CDC)},
  pages 3548--3553, Dec 2016.

\bibitem{FPR14}
M.~Fornasier, B.~Piccoli, and F.~Rossi.
\newblock Mean-field sparse optimal control.
\newblock {\em Phil. Trans. R. Soc. A}, 372(2028):20130400, 2014.

\bibitem{FS}
M.~Fornasier and F.~Solombrino.
\newblock Mean-field optimal control.
\newblock {\em ESAIM: Control, Optimisation and Calculus of Variations},
  20(4):1123--1152, 2014.

\bibitem{FR13}
S.~Foucart and H.~Rauhut.
\newblock {\em A mathematical introduction to compressive sensing}.
\newblock Applied and Numerical Harmonic Analysis. Birkh\"{a}user/Springer, New
  York, 2013.

\bibitem{Fuller}
A.~Fuller.
\newblock Relay control systems optimized for various performance criteria.
\newblock {\em IFAC Proceedings Volumes}, 1(1):520--529, 1960.
\newblock 1st International IFAC Congress on Automatic and Remote Control,
  Moscow, USSR, 1960.

\bibitem{GG2011}
M.~Garavello and P.~Goatin.
\newblock The {A}w-{R}ascle traffic model with locally constrained flow.
\newblock {\em J. Math. Anal. Appl.}, 378(2):634--648, 2011.

\bibitem{GGLP2019}
M.~Garavello, P.~Goatin, T.~Liard, and B.~Piccoli.
\newblock A multiscale model for traffic regulation via autonomous vehicles.
\newblock {\em Journal of Differential Equations}, 296:6088--6124, 2020.

\bibitem{GHPbook}
M.~Garavello, K.~Han, and B.~Piccoli.
\newblock {\em Models for vehicular traffic on networks}, volume~9 of {\em AIMS
  Series on Applied Mathematics}.
\newblock American Institute of Mathematical Sciences (AIMS), Springfield, MO,
  2016.

\bibitem{GPbook}
M.~Garavello and B.~Piccoli.
\newblock {\em Traffic flow on networks}, volume~1 of {\em AIMS Series on
  Applied Mathematics}.
\newblock American Institute of Mathematical Sciences (AIMS), Springfield, MO,
  2006.
\newblock Conservation laws models.

\bibitem{Gasser2013}
I.~Gasser, C.~Lattanzio, and A.~Maurizi.
\newblock Vehicular traffic flow dynamics on a bus route.
\newblock {\em Multiscale Model. Simul.}, 11(3):925--942, 2013.

\bibitem{Gazis1961}
D.~C. Gazis, R.~Herman, and R.~W. Rothery.
\newblock Nonlinear follow-the-leader models of traffic flow.
\newblock {\em Operations research}, 9(4):545--567, 1961.

\bibitem{gloudemans2020interstate}
D.~Gloudemans, W.~Barbour, N.~Gloudemans, M.~Neuendorf, B.~Freeze, S.~ElSaid,
  and D.~B. Work.
\newblock Interstate-24 motion: Closing the loop on smart mobility.
\newblock In {\em 2020 IEEE Workshop on Design Automation for CPS and IoT
  (DESTION)}, pages 49--55. IEEE, 2020.

\bibitem{Golse16}
F.~Golse.
\newblock On the dynamics of large particle systems in the mean field limit.
\newblock In {\em Macroscopic and large scale phenomena: coarse graining, mean
  field limits and ergodicity}, volume~3 of {\em Lect. Notes Appl. Math.
  Mech.}, pages 1--144. Springer, [Cham], 2016.

\bibitem{GomezSerrano2012}
J.~G\'{o}mez-Serrano, C.~Graham, and J.-Y. Le~Boudec.
\newblock The bounded confidence model of opinion dynamics.
\newblock {\em Math. Models Methods Appl. Sci.}, 22(2):1150007, 46, 2012.

\bibitem{GHPV23}
X.~Gong, M.~Herty, B.~Piccoli, and G.~Visconti.
\newblock Crowd dynamics: Modeling and control of multiagent systems.
\newblock {\em Annual Review of Control, Robotics, and Autonomous Systems},
  6(1):null, 2023.

\bibitem{Gong2022}
X.~Gong and A.~Keimer.
\newblock On the well-posedness of the "bando-follow the leader" car following
  model and a time-delayed version".
\newblock {\em {Preprint}}, 2022.
\newblock {R}esearchgate DOI: 10.13140/RG.2.2.22507.62246.

\bibitem{GPV21}
X.~Gong, B.~Piccoli, and G.~Visconti.
\newblock Mean-field of optimal control problems for hybrid model of multilane
  traffic.
\newblock {\em IEEE Control Syst. Lett.}, 5(6):1964--1969, 2021.

\bibitem{greenshields}
B.~{\sc Greenshields}.
\newblock A study of traffic capacity.
\newblock {\em Proceedings of the Highway Research Board}, 14(1):448--477,
  1935.

\bibitem{GUERIAU2016266}
M.~Guériau, R.~Billot, N.-E.~E. Faouzi], J.~Monteil, F.~Armetta, and
  S.~Hassas.
\newblock How to assess the benefits of connected vehicles? a simulation
  framework for the design of cooperative traffic management strategies.
\newblock {\em Transportation Research Part C: Emerging Technologies}, 67:266
  -- 279, 2016.

\bibitem{HHK10}
S.~Y. Ha, T.~Ha, and J.~H. Kim.
\newblock Emergent behavior of a cucker-smale type particle model with
  nonlinear velocity couplings.
\newblock {\em IEEE Transactions on Automatic Control}, 55(7):1679--1683, July
  2010.

\bibitem{Hegselmann2002}
R.~Hegselmann and U.~Krause.
\newblock Opinion dynamics and bounded confidence models, analysis, and
  simulation.
\newblock {\em Journal of Artificial Societies and Social Simulation}, 5(3),
  2002.

\bibitem{Helbing1995}
D.~Helbing and P.~Moln{\'a}r.
\newblock Social force model for pedestrian dynamics.
\newblock {\em Physical Review E}, 51(5):4282--4286, 1995.

\bibitem{Herty2010}
M.~Herty and L.~Pareschi.
\newblock Fokker-{P}lanck asymptotics for traffic flow models.
\newblock {\em Kinet. Relat. Models}, 3(1):165--179, 2010.

\bibitem{Herty2011}
M.~Herty and C.~Ringhofer.
\newblock Averaged kinetic models for flows on unstructured networks.
\newblock {\em Kinet. Relat. Models}, 4(4):1081--1096, 2011.

\bibitem{HR95}
H.~Holden and N.~H. Risebro.
\newblock A mathematical model of traffic flow on a network of unidirectional
  roads.
\newblock {\em SIAM J. Math. Anal.}, 26(4):999--1017, 1995.

\bibitem{Jabin2014}
P.~Jabin and S.~Motsch.
\newblock Clustering and asymptotic behavior in opinion formation.
\newblock {\em Journal of Differential Equations}, 257(11):4165--4187, 12 2014.

\bibitem{jin2013multi}
W.-L. Jin.
\newblock A multi-commodity lighthill-whitham-richards model of lane-changing
  traffic flow.
\newblock {\em Procedia-Social and Behavioral Sciences}, 80:658--677, 2013.

\bibitem{MOBIL07}
A.~Kesting, M.~Treiber, and D.~Helbing.
\newblock General lane-changing model mobil for car-following models.
\newblock {\em Transportation Research Record}, 1999(1):86--94, 2007.

\bibitem{kupka}
I.~A.~K. Kupka.
\newblock The ubiquity of {F}uller's phenomenon.
\newblock In {\em Nonlinear controllability and optimal control}, volume 133 of
  {\em Monogr. Textbooks Pure Appl. Math.}, pages 313--350. Dekker, New York,
  1990.

\bibitem{Lattanzio2011}
C.~Lattanzio, A.~Maurizi, and B.~Piccoli.
\newblock Moving bottlenecks in car traffic flow: a {PDE}-{ODE} coupled model.
\newblock {\em SIAM J. Math. Anal.}, 43(1):50--67, 2011.

\bibitem{laval2006lane}
J.~A. Laval and C.~F. Daganzo.
\newblock Lane-changing in traffic streams.
\newblock {\em Transportation Research Part B: Methodological}, 40(3):251--264,
  2006.

\bibitem{LLG98}
J.-P. Lebacque, J.~B. Lesort, and F.~Giorgi.
\newblock Introducing buses into first-order macroscopic traffic flow models.
\newblock {\em Transportation Reasearch Record}, 1644:70--79, 1998.

\bibitem{Lewin1951}
K.~Lewin.
\newblock {\em Field Theory in Social Science}.
\newblock Harper \& Brothers, 1951.

\bibitem{LiardPiccoli19}
T.~Liard and B.~Piccoli.
\newblock Well-posedness for scalar conservation laws with moving flux
  constraints.
\newblock {\em SIAM J. Appl. Math.}, 79(2):641--667, 2019.

\bibitem{LiardPiccoli21}
T.~Liard and B.~Piccoli.
\newblock On entropic solutions to conservation laws coupled with moving
  bottlenecks.
\newblock {\em Commun. Math. Sci.}, 19(4):919--945, 2021.

\bibitem{Mallat09}
S.~Mallat.
\newblock {\em A wavelet tour of signal processing}.
\newblock Elsevier/Academic Press, Amsterdam, third edition, 2009.
\newblock The sparse way, With contributions from Gabriel Peyr\'{e}.

\bibitem{Maury2011}
B.~Maury and J.~Venel.
\newblock A discrete contact model for crowd motion.
\newblock {\em ESAIM: Mathematical Modelling and Numerical Analysis},
  45(1):145--168, 2011.

\bibitem{Med14}
G.~S. Medvedev.
\newblock The nonlinear heat equation on dense graphs and graph limits.
\newblock {\em SIAM Journal on Mathematical Analysis}, 46(4):2743--2766, 2014.

\bibitem{Motsch2014}
S.~Motsch and E.~Tadmor.
\newblock Heterophilious dynamics enhances consensus.
\newblock {\em SIAM Review}, 56(4):577--621, 2014.

\bibitem{PT22}
T.~Paul and E.~Trélat.
\newblock From microscopic to macroscopic scale equations: mean field,
  hydrodynamic and graph limits, 2022.

\bibitem{Payne1971}
H.~{\sc Payne}.
\newblock Models of freeway traffic and control.
\newblock {\em Mathematical models of public systems}, 1(1):51--61, 1971.

\bibitem{Perthame04}
B.~Perthame.
\newblock Mathematical tools for kinetic equations.
\newblock {\em Bull. Amer. Math. Soc. (N.S.)}, 41(2):205--244, 2004.

\bibitem{PC19book}
G.~Peyré and M.~Cuturi.
\newblock Computational optimal transport: With applications to data science.
\newblock {\em Foundations and Trends® in Machine Learning}, 11(5-6):355--607,
  2019.

\bibitem{PiacentiniGoatinFerrara2019platoon}
G.~Piacentini, P.~Goatin, and A.~Ferrara.
\newblock {A macroscopic model for platooning in highway traffic}.
\newblock working paper or preprint, Oct. 2019.

\bibitem{piccoli2019}
B.~Piccoli.
\newblock Measure differential equations.
\newblock {\em Arch. Ration. Mech. Anal.}, 233(3):1289--1317, 2019.

\bibitem{PPT19}
B.~Piccoli, N.~P. Duteil, and E.~Tr\'{e}lat.
\newblock Sparse control of hegselmann--krause models: Black hole and
  declustering.
\newblock {\em SIAM Journal on Control and Optimization}, 57(4):2628--2659,
  2019.

\bibitem{PR13}
B.~Piccoli and F.~Rossi.
\newblock Transport equation with nonlocal velocity in {W}asserstein spaces:
  convergence of numerical schemes.
\newblock {\em Acta Appl. Math.}, 124:73--105, 2013.

\bibitem{PR14}
B.~Piccoli and F.~Rossi.
\newblock Generalized {W}asserstein distance and its application to transport
  equations with source.
\newblock {\em Arch. Ration. Mech. Anal.}, 211(1):335--358, 2014.

\bibitem{PR16}
B.~Piccoli and F.~Rossi.
\newblock On properties of the generalized {W}asserstein distance.
\newblock {\em Arch. Ration. Mech. Anal.}, 222(3):1339--1365, 2016.

\bibitem{PR21}
B.~Piccoli and F.~Rossi.
\newblock Generalized solutions to bounded-confidence models.
\newblock {\em Mathematical Models and Methods in Applied Sciences},
  31(06):1237--1276, 2021.

\bibitem{PT09}
B.~Piccoli and A.~Tosin.
\newblock Pedestrian flows in bounded domains with obstacles.
\newblock {\em Contin. Mech. Thermodyn.}, 21(2):85--107, 2009.

\bibitem{PT11}
B.~Piccoli and A.~Tosin.
\newblock Time-evolving measures and macroscopic modeling of pedestrian flow.
\newblock {\em Archive for Rational Mechanics and Analysis}, 199(3):707--738,
  2011.

\bibitem{PT12}
B.~Piccoli and A.~Tosin.
\newblock Vehicular traffic: a review of continuum mathematical models.
\newblock In {\em Mathematics of complexity and dynamical systems. {V}ols.
  1--3}, pages 1748--1770. Springer, New York, 2012.

\bibitem{RamadanSeibold2017}
R.~A. Ramadan and B.~Seibold.
\newblock Traffic flow control and fuel consumption reduction via moving
  bottlenecks.
\newblock In {\em Transportation Research Board Conference}, 2017.

\bibitem{STERN2019351}
R.~E. Stern, Y.~Chen, M.~Churchill, F.~Wu, M.~L.~D. Monache], B.~Piccoli,
  B.~Seibold, J.~Sprinkle, and D.~B. Work.
\newblock Quantifying air quality benefits resulting from few autonomous
  vehicles stabilizing traffic.
\newblock {\em Transportation Research Part D: Transport and Environment},
  67:351 -- 365, 2019.

\bibitem{stern2018dissipation}
R.~E. Stern, S.~Cui, M.~L. Delle~Monache, R.~Bhadani, M.~Bunting, M.~Churchill,
  N.~Hamilton, H.~Pohlmann, F.~Wu, B.~Piccoli, et~al.
\newblock Dissipation of stop-and-go waves via control of autonomous vehicles:
  Field experiments.
\newblock {\em Transp. Research Part C: Emerging Techn.}, 89:205--221, 2018.

\bibitem{sugiyama2008}
Y.~Sugiyama, M.~Fukui, M.~Kikuchi, K.~Hasebe, A.~Nakayama, K.~Nishinari, S.-i.
  Tadaki, and S.~Yukawa.
\newblock Traffic jams without bottlenecks—experimental evidence for the
  physical mechanism of the formation of a jam.
\newblock {\em New journal of physics}, 10(3):033001, 2008.

\bibitem{TalebBS}
N.~Taleb.
\newblock {\em The Black Swan - the impact of the highly improbable}.
\newblock Random House, New York, 2007.

\bibitem{TALEBPOUR2016143}
A.~Talebpour and H.~S. Mahmassani.
\newblock Influence of connected and autonomous vehicles on traffic flow
  stability and throughput.
\newblock {\em Transportation Research Part C: Emerging Technologies}, 71:143
  -- 163, 2016.

\bibitem{MR2247927}
G.~Toscani.
\newblock Kinetic models of opinion formation.
\newblock {\em Commun. Math. Sci.}, 4(3):481--496, 2006.

\bibitem{Treiber2000}
M.~Treiber, A.~Hennecke, and D.~Helbing.
\newblock Congested traffic states in empirical observations and microscopic
  simulations.
\newblock {\em Physical review E}, 62(2):1805, 2000.

\bibitem{treiterer1974}
J.~Treiterer and J.~Myers.
\newblock The hysteresis phenomenon in traffic flow.
\newblock {\em Transportation and traffic theory}, 6:13--38, 1974.

\bibitem{villani}
C.~Villani.
\newblock {\em Topics in optimal transportation}, volume~58 of {\em Graduate
  Studies in Mathematics}.
\newblock American Mathematical Society, Providence, RI, 2003.

\bibitem{vi09}
C.~Villani.
\newblock {\em Optimal Transport}, volume 338 of {\em Grundlehren der
  Mathematischen Wissenschaften [Fundamental Principles of Mathematical
  Sciences]}.
\newblock Springer-Verlag, Berlin, 2009.
\newblock Old and new.

\bibitem{WW11}
G.~Wachsmuth and D.~Wachsmuth.
\newblock Convergence and regularization results for optimal control problems
  with sparsity functional.
\newblock {\em ESAIM Control Optim. Calc. Var.}, 17(3):858--886, 2011.

\bibitem{7362183}
M.~{Wang}, W.~{Daamen}, S.~P. {Hoogendoorn}, and B.~{van Arem}.
\newblock Cooperative car-following control: Distributed algorithm and impact
  on moving jam features.
\newblock {\em IEEE Transactions on Intelligent Transportation Systems},
  17(5):1459--1471, 2016.

\bibitem{Whitham74}
G.~B. Whitham.
\newblock {\em Linear and nonlinear waves}.
\newblock Wiley-Interscience, New York, 1974.

\bibitem{WU201982}
F.~Wu, R.~E. Stern, S.~Cui, M.~L.~D. Monache], R.~Bhadani, M.~Bunting,
  M.~Churchill, N.~Hamilton, R.~Haulcy, B.~Piccoli, B.~Seibold, J.~Sprinkle,
  and D.~B. Work.
\newblock Tracking vehicle trajectories and fuel rates in phantom traffic jams:
  Methodology and data.
\newblock {\em Transportation Research Part C: Emerging Technologies}, 99:82 --
  109, 2019.

\bibitem{YKLH17}
K.~Yuan, V.~L. Knoop, L.~Leclercq, and S.~P. Hoogendoorn.
\newblock Capacity drop: a comparison between stop-and-go wave and standing
  queue at lane-drop bottleneck.
\newblock {\em Transportmetrica B: Transport Dynamics}, 5(2):145--158, 2017.

\bibitem{ZHENG201416}
Z.~Zheng.
\newblock Recent developments and research needs in modeling lane changing.
\newblock {\em Transportation Research Part B: Methodological}, 60:16--32,
  2014.

\end{thebibliography}


\end{document}